\documentclass{article}

\usepackage[T1]{fontenc}
\usepackage{verbatim}
\usepackage{amsthm, amsmath}
\usepackage{indentfirst}
\usepackage{graphicx}
\usepackage{stfloats}
\usepackage{enumerate}
\usepackage[margin=2.5cm]{geometry}
\usepackage[]{caption2}
\usepackage{hyperref}

\newtheorem{lemma2026ErdosQuestion}{Lemma}[section]
\newtheorem{theorem2026ErdosQuestion}[lemma2026ErdosQuestion]{Theorem}
\newtheorem{conjecture2026ErdosQuestion}[lemma2026ErdosQuestion]{Conjecture}
\newtheorem{definition2026ErdosQuestion}[lemma2026ErdosQuestion]{Definition}
\newtheorem{question2026ErdosQuestion}[lemma2026ErdosQuestion]{Question}
\newtheorem{claim2026ErdosQuestion}{Claim}

\begin{document}
\title{An extremal question of Erd\H{o}s et al. on cycles}\date{}
\author{
 Zhiyong Gan$^{1}$\thanks{E-mail address: ganzhiyong@m.scnu.edu.cn},\quad
 Zhiliang Tan$^{1}$,\quad
 Yi Huang$^{2}$\thanks{Corresponding author. E-mail address: yihuang@m.scnu.edu.cn}
\\\\1 School of Computer Science, South China Normal University, \\Guangzhou 510631, P. R. China
\\\\2 Aberdeen Institute of Data Science and Artificial Intelligence, \\South China Normal University, \\Guangzhou 510631, P. R. China
}
\maketitle
\begin{abstract}
This paper addresses an extremal problem on long cycles proposed by Erd\H{o}s, Faudree, Schelp and Simonovits in 1989. Let positive integers satisfy $2\leq k\leq d<n$. For an $n$-vertex simple graph $G$, denote by $\theta_+^G$ the number of vertices of degree at least $d$, and let $C_{\geq k+1}$ stand for a cycle of length at least $k+1$. The original question conjectured that $G$ contains a $C_{\geq k+1}$ whenever
$$
    \theta_+^G \geq\left\lfloor\frac{k}{2}\right\rfloor \cdot\left\lfloor\frac{n-1}{d}\right\rfloor+2.
$$
However, Li, Ma and Ning disproved this conjecture by constructing infinite families of counterexamples and provided linear lower bounds for reference. Later, Gan, Yi and Jiang established piecewise sufficient conditions depending on the parity of $k$, yet their bounds are sharp only in the case $k=d$ and remain loose for $k\neq d$.

In this work, we derive a sufficient condition for the existence of $C_{\geq k+1}$. For all admissible parameters except the special case that $k$ is odd and $k=d$, our threshold formula is unified and does not depend on the parity of $k$. Our threshold has two branches: the exceptional case $k=d$ odd, where we adopt Woodall's sharp threshold, and the unified formula for all other admissible parameters, which does not depend on the parity of $k$. Our main theorem states the following: If
$$
    \theta_+^G\geq
    \begin{cases}
        \displaystyle \frac{k(n-2)}{2(k-1)}+1 & \text{$k$ odd and } k=d,\\[6pt]
        \displaystyle\left\lfloor \frac {k}{2}\right\rfloor \cdot \frac {n-1}{d}
        +\frac{(n-d-1)\lfloor k/2\rfloor}{d(d-\lfloor k/2\rfloor)}+1 & \text{otherwise},
    \end{cases}
$$
then $G$ contains a cycle $C_{\geq k+1}$. Moreover, we construct a family of $C_{\geq k+1}$-free graphs to verify that the above bound is tight for all admissible integers $k$. This paper solves the long-standing Erd\H{o}s-type extremal problem on high-degree vertices and long cycles.


\end{abstract}

\section{Introduction}
\subsection{Background and Classical Extremal Results}
\noindent A central research thread in extremal graph theory concerns the existence of long paths and cycles constrained by vertex degree distributions. A core fundamental question reads: given an $n$-vertex simple undirected graph, how many vertices of degree at least $d$ are required to force a cycle of length at least $k+1$ (denoted $C_{\geq k+1}$)?

Two landmark classical results lay the foundation for this direction. First, Dirac's cycle theorem \cite{articleDirac1952} asserts that a graph with minimum degree at least $d$ contains a cycle $C_{\geq d+1}$. Second, the Erd\H{o}s-Gallai theorem \cite{articleErdosGallai1959} provides an edge-count threshold ensuring the presence of a cycle $C_{\geq k+1}$. Both results characterize cycle existence via global minimum degree or total edge number, but they fail to handle general graphs where only a subset of vertices attains high degree.

In 1975, Woodall \cite{articleWoodall1975} pioneered the study of the counting parameter $\theta_+^G$, the total number of vertices with degree no less than $d$, and derived sufficient criteria for long cycles solely based on $\theta_+^G$. His work became the primary tool for subsequent Erd\H{o}s-type conjectures. Later, Jackson \cite{articleJackson1981} followed by Kostochka, Luo and Zirlin \cite{articleKostochka2020}, proved pivotal lemmas on bipartite and 2-connected bipartite graphs that guarantee cycles covering specified vertex sets; these tools are indispensable for analyzing the neighborhood structure of high-degree vertices.

\subsection{The Erd\H{o}s Open Questions and Prior Literature}
At the 1986 Jinan Conference on Graph Theory and Its Applications (published in \textit{Annals of the New York Academy of Sciences}, 1989, \cite{articleErdosEtal1989}), Erd\H{o}s, Faudree, Schelp and Simonovits proposed three related extremal problems: two for long paths $P_{k+1}$ and one for long cycles $C_{\geq k+1}$. This paper mainly solves the cycle-related open problem, stated as Question 1.2 below. We also briefly introduce the path conjecture and its counterexamples for completeness.

\begin{conjecture2026ErdosQuestion}
[Erd\H{o}s et al., 1989 \cite{articleErdosEtal1989}]\label{conjectureErdosEtal1989}
    Let $n,d,k$ be positive integers with $k\leq d<n$, and let $G$ be an $n$-vertex graph. Let
    $$
        t=
            \begin{cases}
                1 & k \text{ odd},\\
                2 & k \text{ even}.
            \end{cases}
    $$
    If
    $$
        \theta_+^G\geq \left\lfloor\frac{k-1}{2}\right\rfloor\cdot \left\lfloor\frac{n}{d+1}\right\rfloor + t,
    $$
    then $G$ contains a path $P_{k+1}$.
\end{conjecture2026ErdosQuestion}
Li, Ma and Ning (2022) \cite{articleLiMaNing2022} proved that Conjecture~\ref{conjectureErdosEtal1989} holds for all integers $k\ne 4$ and fails when $k=4$.

\begin{question2026ErdosQuestion}
[Erd\H{o}s et al., 1989, \cite{articleErdosEtal1989}]\label{questionsErdosEtal1989A}
    Let $n,d,k$ be positive integers with $2\leq k\leq d<n$, and let $G$ be an $n$-vertex graph. Suppose
    $$
        \theta_+^G\geq \left\lfloor\frac{k}{2}\right\rfloor\cdot \left\lfloor\frac{n-1}{d}\right\rfloor + 2.
    $$
    Does it follow that $G$ contains a cycle $C_{\geq k+1}$?
\end{question2026ErdosQuestion}
The authors proposed a construction intuition: glue $\left\lfloor (n-1)/d \right\rfloor$ copies of the graph $H_{d,k+1}$ by identifying one vertex per copy. Each copy is $C_{\geq k+1}$-free and contains exactly $\left\lfloor k/2 \right\rfloor$ high-degree vertices.

Subsequent studies exposed the insufficiency of this conjectured linear bound:
\begin{enumerate}[(1)]
\item Li, Ma and Ning (2022) \cite{articleLiMaNing2022} constructed infinite families of counterexamples to Question \ref{questionsErdosEtal1989A}. They showed for any fixed $k\geq 2$, there exists a constant $c>0$ such that infinitely many $C_{\geq k+1}$-free graphs satisfy
    $$
        \theta_+^G>\left(\left\lfloor\frac{k}{2}\right\rfloor +c\right)\cdot \frac{n-1}{d}.
    $$
\item Woodall proved two foundational theorems for long cycles with high-degree vertex counts:
\begin{theorem2026ErdosQuestion}
[Woodall, 1975, \cite{articleWoodall1975}]\label{theoremWoodall1975NDD}
    Let $n,d$ be positive integers with $2 \leq d < n$ and $G$ be an $n$-vertex graph. If
    $$
        \theta_+^G\ge
        \begin{cases}
            \dfrac{(d+2)(n-1)}{2d} & d\text{ even},\\[4pt]
            \dfrac{d(n-2)}{2(d-1)}+1 & d\text{ odd},
        \end{cases}
    $$
    then $G$ contains a cycle $C_{\geq d+1}$.
\end{theorem2026ErdosQuestion}

\begin{theorem2026ErdosQuestion}
[Woodall, 1975, \cite{articleWoodall1975}]\label{theoremWoodall1975NKD}
    Let $n,k,d$ be positive integers with $2\leq k\leq d<n$ and $G$ be an $n$-vertex graph. If
    $$
        \theta_+^G\geq \frac{(k+3)(n-1)}{2d},
    $$
    then $G$ contains a cycle $C_{\geq k+1}$.
\end{theorem2026ErdosQuestion}

\item Gan, Yi and Jiang generalized Woodall's results with improved piecewise bounds:
\begin{theorem2026ErdosQuestion}
[Gan et al., submitted, \cite{articleGanYiJiang2026}]\label{theoremMyResult1}
    Let $n,k,d$ be positive integers with $2\leq k\leq d<n$ and $G$ be an $n$-vertex graph. If
    $$
        \theta_+^G\ge
        \begin{cases}
            \dfrac{(k+2)(n-1)}{2d} & k\text{ even},\\[4pt]
            \dfrac{k(n-2)}{2(d-1)}+1 & k\text{ odd},
        \end{cases}
    $$
    then $G$ contains a $C_{\geq k+1}$ cycle.
\end{theorem2026ErdosQuestion}
This bound recovers the two branches of Theorem~\ref{theoremWoodall1975NDD} as the special case $k=d$, but is non-sharp whenever $k\ne d$.
\end{enumerate}

Erd\H{o}s et al. also posed a connected-graph version for long paths:
\begin{question2026ErdosQuestion}
[Erd\H{o}s et al., 1989, \cite{articleErdosEtal1989}]\label{questionsErdosEtal1989B}
    Let $n, d, k$ be positive integers with $k \leq d < n$ and $G$ be a connected $n$-vertex graph. Suppose
    $$
        \theta_+^G\geq \left\lfloor\frac{k-3}{4}\right\rfloor\cdot \left\lfloor\frac{n-1}{d}\right\rfloor +2.
    $$
    Does $G$ contain a path $P_{k + 1}$?
\end{question2026ErdosQuestion}
Li et al. disproved Question \ref{questionsErdosEtal1989B} for all $k\ge7$ via glued star-template counterexamples. Inspired by their construction, we propose the following conjecture for connected graphs:
\begin{conjecture2026ErdosQuestion}
\label{conjectureOfMyErdosQuestion}
    Let $n, d, k$ be positive integers with $k \leq d < n$ and $G$ be a connected $n$-vertex graph. If
    $$
        \theta_+^G\geq \left\lfloor\frac{k-3}{4}\right\rfloor\cdot \frac{n-1}{d}
        +\left(1-\frac1d\left\lfloor\frac{k-3}{4}\right\rfloor\right)
        \frac{(n-1)\lfloor (k-3)/4\rfloor}
        {\lfloor (k-3)/4\rfloor+d(d-1)} +2,
    $$
    then $G$ contains a path $P_{k + 1}$.
\end{conjecture2026ErdosQuestion}

\subsection{Main Result and Paper Organization}
All prior sufficient bounds are piecewise and split by the parity of $k$; furthermore, those bounds are tight only in the special case $k=d$. Our core contribution is a new tight threshold. For all parameters with $2\leq k\leq d<n$ except the case when $k$ is odd and $k=d$, our formula is unified and free of parity distinction on $k$. For the exceptional case where $k$ is odd and $k=d$, we adopt the known sharp threshold from Woodall. The full main theorem is restated below.
\begin{theorem2026ErdosQuestion}[Main Result]\label{theorem:main}
    Let $n,k,d$ be positive integers with $2\leq k\leq d<n$ and $G$ be an $n$-vertex graph. If
    $$
        \theta_+^G\geq
        \begin{cases}
            \displaystyle \frac{k(n-2)}{2(k-1)}+1 & k\text{ odd and } k=d,\\[6pt]
            \displaystyle\left\lfloor \frac {k}{2}\right\rfloor \cdot \frac {n-1}{d}
            +\frac{(n-d-1)\lfloor k/2\rfloor}{d(d-\lfloor k/2\rfloor)}+1 & \text{otherwise},
        \end{cases}
    $$
    then $G$ contains a cycle $C_{\geq k+1}$.
\end{theorem2026ErdosQuestion}
We construct an explicit family of $C_{\geq k+1}$-free graphs achieving equality in the inequality of Theorem~\ref{theorem:main}, which certifies tightness.
\begin{itemize}
    \item If $d$ is odd and $k=d$, then we define the examples $O(d, t)$ recursively: Let $O(d, 1)$ be a graph on $2 d$ vertices consisting of $d$ vertices of degree $d$ forming a $K_d$, each of which is joined to one of $d$ vertices of degree 1; $O(d, t)$ is obtained from $O(d, 1)$ and $O(d, t-1)$ by identifying a terminal vertex of $O(d, t-1)$, and its adjacent vertex of degree $d$, with a similar pair of vertices taken from $O(d, 1)$, but in the reverse order. (So $O(d, t)$ has $t(2 d-2)+2$ vertices, all of degree $d$ or 1, and $\frac{k(n-2)}{2(k-1)}$ high-degree vertices.) This family of $C_{\geq k+1}$-free graphs was first introduced by Woodall \cite{articleWoodall1975}.

    \item Otherwise, let $\alpha, \beta$ be positive integers satisfying $d=(1+\alpha)\left\lfloor\frac{k}{2}\right\rfloor$ and $n=1+d+\alpha \beta d$; in particular $2 \leq k \leq d < n$. Let $H_0$ be a copy of $H_{d, k+1}$ with a low-degree vertex $v_0$. For each $1 \leq i \leq \beta$, let $H_i$ be obtained from $\alpha$ copies of $H_{d, k+1}$ by identifying one low-degree vertex from each copy of $H_{d, k+1}$ (call the resulting vertex $u_i$ ); let $v_i$ be a low-degree vertex in $H_i$ other than $u_i$. Finally, let graph $G$ be obtained from $H_0, H_1, \ldots, H_\beta$ by identifying $v_{i-1}$ and $u_i$ for all $1 \leq i \leq \beta$. Since each block of $G$ is a copy of $H_{d, k+1}$, we see that $G$ is $C_{\geq k+1}$-free. However, $G$ has
        $$
            (1 + \alpha \beta)\cdot \left\lfloor\frac{k}{2}\right\rfloor + \beta
            = \left\lfloor\frac{k}{2}\right\rfloor\cdot \frac{n - 1}{d}
            + \frac{(n - d - 1)\left\lfloor k/2\right\rfloor}{d(d-\left\lfloor k/2\right\rfloor)}
        $$
        high-degree vertices. This family of $C_{\geq k+1}$-free graphs was first introduced by Li et al. \cite{articleLiMaNing2022}.

\end{itemize}

The rest of this paper is structured as follows:
\begin{itemize}
    \item Section \ref{sectionDefinition}: Standard graph-theoretic definitions and notation.
    \item Section \ref{sectionPreliminary}: Preliminary lemmas from Dirac, Jackson, Kostochka--Luo--Zirlin used in the proof.
    \item Section \ref{sectionProof}: Full proof of Theorem~\ref{theorem:main} via minimal counterexample analysis, split into two exhaustive structural cases for high-degree vertices.
    \item References: Cited literature.
\end{itemize}

\section{Definitions and Notation}\label{sectionDefinition}
Let $G$ be a finite, simple, undirected graph. Write $V(G)$ for the vertex set, $E(G)$ for the edge set.
\begin{definition2026ErdosQuestion}
Given $v\in V(G)$, let $\deg_G(v)$ denote the \emph{degree} of vertex $v$ in $G$, i.e., the number of edges of $G$ incident to $v$. A vertex $v\in V(G)$ is called \textit{high-degree} if $\deg_G(v)\geq d$, otherwise \textit{low-degree}. Denote:
\begin{itemize}
    \item $D_+^G = \{v\in V(G)\mid \deg_G(v)\geq d\}$ (set of high-degree vertices),
    \item $D_-^G = \{v\in V(G)\mid \deg_G(v)<d\}$ (set of low-degree vertices),
    \item $\theta_+^G = |D_+^G|$.
\end{itemize}
For $v\in V(G)$, $N_G(v)=\{u\in V(G)\mid uv\in E(G)\}$ is the open neighborhood of $v$.
\end{definition2026ErdosQuestion}
\begin{definition2026ErdosQuestion}
$C_m$ denotes a cycle of exactly $m$ vertices; $C_{\geq t}$ denotes any cycle with at least $t$ vertices. A graph is \textit{$C_{\geq t}$-free} if it contains no $C_{\geq t}$. Similarly, $P_m$ is a path on $m$ vertices.%, and $P_{\geq t}$-free graphs contain no such path.
\end{definition2026ErdosQuestion}
Given two vertices $s,t$ of a graph $G$, an \emph{s-t path} in $G$ is a simple path in $G$ with endpoints $s$ and $t$.

Given a cycle $C$ with cyclic ordering $a_1,a_2,\dots,a_m$, write $\overrightarrow{C}$ for the clockwise orientation and $\overleftarrow{C}$ for the reverse. The segment $a_i\overrightarrow{C}a_j$ means the consecutive vertices from $a_i$ to $a_j$ along $\overrightarrow{C}$; $a_j\overleftarrow{C}a_i$ is the reversed segment. The same notation applies to paths.

Standard graph families:
\begin{itemize}
    \item $K_s$: complete graph on $s$ vertices;
    \item $\overline{K}_s$: independent set of $s$ vertices;
    \item $K_{1,s}$: star graph with $s$ leaves.
\end{itemize}
\begin{definition2026ErdosQuestion}
The template graph $H_{d,k}$ is formed by taking disjoint $K_{\left\lfloor (k-1)/2 \right\rfloor}$ and $\overline{K}_{d+1-\left\lfloor (k-1)/2 \right\rfloor}$, then joining every vertex of the clique to every vertex of the independent set.
\end{definition2026ErdosQuestion}
A \textit{block} of $G$ is a maximal 2-connected subgraph of $G$. For any real $x$, $\left\lfloor x\right\rfloor$ denotes the floor of $x$, and $\left\lceil x\right\rceil$ denotes the ceiling of $x$. For brevity, in inline formulas we write $\lfloor k/2\rfloor$ and $\lfloor (k-3)/4\rfloor$ for $\left\lfloor\frac{k}{2}\right\rfloor$ and $\left\lfloor\frac{k-3}{4}\right\rfloor$, respectively. All undefined terminology follows Bondy and Murty's textbook \textit{Graph Theory} \cite{bookBondy2008}.

\section{Preliminary Results} \label{sectionPreliminary}
We collect classical lemmas critical for proving Theorem~\ref{theorem:main}.

\begin{lemma2026ErdosQuestion}
[Dirac, 1952, \cite{articleDirac1952}]\label{lemmaDirac}
    Let $d>1$. If every vertex of $G$ has degree at least $d$, then $G$ contains a cycle $C_{\geq d+1}$.
\end{lemma2026ErdosQuestion}

\begin{lemma2026ErdosQuestion}
[Jackson, 1981, \cite{articleJackson1981}]\label{lemmaJackson1981}
    Let $d$ be a positive integer and $G$ be a bipartite graph with partitions $X, Y$. If $2 \leq |X| \leq d$, $|Y| \leq 2d-2$, and every vertex in $X$ has degree at least $d$, then $G$ has a cycle containing all vertices in $X$.
\end{lemma2026ErdosQuestion}

Lemma~\ref{lemma}, conjectured by Jackson and proved by Kostochka, Luo and Zirlin, strengthens Lemma~\ref{lemmaJackson1981} on 2-connected graphs.
\begin{lemma2026ErdosQuestion}
[Kostochka et al., 2020, \cite{articleKostochka2020}]\label{lemma}
    Let $d$ be a positive integer and $G$ be a 2-connected bipartite graph with partitions $X, Y$. If $2 \leq |X| \leq d$, $|Y| \leq 3d - 5$, and every vertex in $X$ has degree at least $d$, then $G$ has a cycle containing all vertices in $X$.
\end{lemma2026ErdosQuestion}

\section{Proof of Theorem~\ref{theorem:main}}\label{sectionProof}
\begin{proof}
    We first treat the special case $k = d$ independently using Woodall's classic results; after verifying this case holds, we only need to conduct the minimal counterexample argument under the assumption $2 \leq k < d$. If $k$ is even and $k=d$, then
    \begin{align*}
        \theta_+^G &\geq \left\lfloor\frac{k}{2}\right\rfloor\cdot \frac{n - 1}{d}
        + \frac{(n - d - 1)\left\lfloor k/2\right\rfloor}{d(d-\left\lfloor k/2\right\rfloor)} + 1 \\
                    &= \frac{n - 1}{2} + \frac{n - 1}{d} \\
                    &= \frac{(d + 2)(n - 1)}{2d}.
    \end{align*}
    This inequality exactly matches the hypothesis of Woodall's Theorem~\ref{theoremWoodall1975NDD}, which guarantees a cycle $C_{d+1} = C_{k+1}$ in $G$.

    If $k$ is odd and $k=d$, Theorem~\ref{theorem:main} directly adopts the piecewise threshold $$\theta_+^G \geq \dfrac{k(n-2)}{2(k-1)}+1 = \dfrac{d(n-2)}{2(d-1)}+1,$$ coinciding with the odd-$d$ branch of Woodall's Theorem~\ref{theoremWoodall1975NDD}, so $G$ also contains a cycle $C_{d+1} = C_{k+1}$.

    In both subcases with $k=d$, the statement of Theorem~\ref{theorem:main} holds. For the remainder of the proof, we assume $2 \leq k < d$ and proceed by contradiction via minimal counterexample.

    To the contrary, assume for contradiction that there exists an $n$-vertex simple graph $G$ satisfying the following four minimal counterexample conditions:
    \begin{enumerate}[(1)]
        \item The high-degree vertex count satisfies the threshold inequality
            $$
                \theta_+^G\geq \left\lfloor\frac{k}{2}\right\rfloor\cdot \frac{n-1}{d}
                +\frac{(n-d-1)\left\lfloor k/2\right\rfloor}{d(d-\left\lfloor k/2\right\rfloor)}+1,
            $$\label{psiItem1}
        \item $G$ is $C_{\geq k+1}$-free (no cycle of length at least $k+1$);\label{psiItem2}
        \item Subject to Conditions \eqref{psiItem1} and \eqref{psiItem2}, the vertex number $n = |V(G)|$ is minimized among all such graphs; \label{psiItem3}
        \item Subject to Conditions \eqref{psiItem1}, \eqref{psiItem2} and \eqref{psiItem3}, the edge number $|E(G)|$ is minimized among all $n$-vertex graphs.\label{psiItem4}
    \end{enumerate}

    First, minimality implies $G$ has no isolated vertices, and by definition $\theta_+^G = |D_+^G|$. We establish two foundational structural claims on the neighborhood of low-degree vertices.

    \begin{claim2026ErdosQuestion}\label{claimNeighborOfLow}
        The graph $G$ satisfies the following two structural properties:
        \begin{enumerate}[(i)]
            \item $G$ contains at least one low-degree vertex in $D_-^G$; \label{claimNeighborOfLow1}
            \item Every neighbor of any low-degree vertex in $D_-^G$ is a high-degree vertex in $D_+^G$.\label{claimNeighborOfLow2}
        \end{enumerate}
    \end{claim2026ErdosQuestion}
    \begin{proof}
        \emph{(i)} Suppose for contradiction that every vertex in $G$ is high-degree, so $\deg_G(v) \geq d$ for all $v\in V(G)$. By Dirac's cycle lemma (Lemma~\ref{lemmaDirac}), $G$ admits a cycle $C_{\geq d+1}$. Since $k < d$, we have $C_{\geq d+1}\subseteq C_{\geq k+1}$, which contradicts the $C_{\geq k+1}$-free assumption of $G$. Therefore low-degree vertices must exist.

        \emph{(ii)} Let $u\in D_-^G$ be an arbitrary low-degree vertex, and take any neighbor $v\in N_G(u)$. Assume to the contrary that $v\in D_-^G$. Construct a new graph $G' = G - uv$ by deleting the edge $uv$. We verify that $G'$ still satisfies Conditions \eqref{psiItem1} and \eqref{psiItem2}:
        \begin{itemize}
            \item Deleting one incident edge of two low-degree vertices does not change the set $D_-^G$, hence $\theta_+^{G'} = \theta_+^G$ and Condition \eqref{psiItem1} holds;
            \item Removing an edge cannot create a long cycle, so $G'$ remains $C_{\geq k+1}$-free, satisfying Condition \eqref{psiItem2}.
        \end{itemize}
        The graph $G'$ shares the same vertex set $V(G')=V(G)$ but has strictly fewer edges than $G$, violating the edge-minimality Condition \eqref{psiItem4}. We conclude $v\in D_+^G$, completing the proof of Claim~\ref{claimNeighborOfLow}.
    \end{proof}\smallskip

    \begin{claim2026ErdosQuestion}\label{claim:k2}
        Theorem~\ref{theorem:main} holds for $k=2$.
    \end{claim2026ErdosQuestion}

    \begin{proof}
        Let $G$ be a minimal counterexample, i.e., $G$ is $C_{\geq 3}$-free (hence a forest) and satisfies
        \begin{equation}\label{eq:k2-threshold}
        \theta_+^G\geq \frac{n-1}{d}+\frac{n-d-1}{d(d-1)}+1
        =\frac{n-2}{d-1}+1.
        \end{equation}
        Since $d<n$, we have $\frac{n-2}{d-1}\geq 1$, so $\theta_+^G\geq 2$. Hence there exists a high-degree vertex $v\in D_+^G$. Since $G$ is a forest, starting from $v$ and following a maximal path avoiding previously visited vertices, we eventually reach a leaf $w$. By Claim~\ref{claimNeighborOfLow}~\eqref{claimNeighborOfLow2}, the unique neighbor $v'$ of $w$ is a high-degree vertex, so $\deg(v')\geq d$.

        Consider the neighbors of $v'$ other than the one on the path from $v$ to $w$. There are at least $d-1$ such neighbors, and one of them is $w$. If all of them are leaves, then deleting these $d-1$ leaves (including $w$) reduces $n$ by $d-1$ and reduces $\theta_+^G$ by at most $1$ (since $v'$ may become low-degree, but no other high-degree vertex is affected). The resulting graph $G'$ has $n'=n-(d-1)$ vertices and satisfies
        $$
            \theta_+^{G'}\geq\theta_+^G-1.
        $$
        By \eqref{eq:k2-threshold},
        $$
            \theta_+^{G'}\geq \frac{n-2}{d-1} = \frac{n-d-1}{d-1}+1 = \frac{n'-2}{d-1}+1.
        $$
        So $G'$ satisfies the threshold condition. Since $G'$ is still $C_{\geq 3}$-free and has fewer vertices, this contradicts the minimality of $G$.
        
        Hence some neighbor $w_1$ of $v'$ (other than the path neighbor) is not a leaf. Then $w_1$ has another neighbor $x\neq v'$, and we continue the argument with $x$. Since $G$ is finite, this process must terminate at a high-degree vertex $v''$ all of whose non-path neighbors are leaves. Deleting these $d-1$ leaves again yields a smaller counterexample, a contradiction.

        Therefore $G$ must contain a cycle $C_{\geq 3}$, and the theorem holds for $k=2$.
    \end{proof}\smallskip

    By Claim~\ref{claim:k2}, we have $k>2$.

    \begin{claim2026ErdosQuestion}\label{claimLowerN}
        The graph $G$ satisfies the two lower bounds below:
        \begin{enumerate}[(i)]
            \item $\theta_+^G \geq \left\lfloor\frac{k}{2}\right\rfloor + 1$, and\label{claimItem1}
            \item $n \geq \frac{3}{2}d + 1$.\label{claimItem2}
        \end{enumerate}
    \end{claim2026ErdosQuestion}
    \begin{proof}
        \emph{(i)} Since $d < n$, we have $\frac{n-1}{d} \geq 1$. Hence $\theta_+^G \geq \left\lfloor\frac{k}{2}\right\rfloor + 1$, that is, \emph{(i)} holds.

        \emph{(ii)} We prove \emph{(ii)} by contrapositive. Assume instead $n \leq \tfrac{3}{2}d$, and split the argument based on the parity of $k$ (recall $k < d$ throughout this claim).

        \textbf{Case A.} $k$ is even.

            Write $|X| = \frac{k}{2}+1$ for a subset $X\subseteq D_+^G$, and set $Y = V(G)\setminus X$. Let $G' = G(X, Y)$ denote the spanning bipartite graph of $G$ with bipartition $(X,Y)$.

            Define $d' = d-\frac{k}{2}$. For each $x\in X$, its degree within the bipartite subgraph satisfies $$\deg_{G'}(x) \geq \deg_G(x) - (|X| - 1) \geq d - \frac{k}{2} = d'.$$

            Since $2\leq k< d$, we have $2\leq |X|\leq d'$. By our contrary assumption $n \leq \frac{3}{2}d$, $$|Y| = n - |X| \leq \frac{3}{2}d-(\frac{k}{2} + 1) = 2d - k - 2 + \frac{k - d + 2}{2} \leq 2d - k - 2 + \frac{1}{2}.$$ Since $|Y|$ is an integer, $|Y| \leq 2d-k-2 = 2d'-2$. By Jackson's Lemma~\ref{lemmaJackson1981}, $G$ contains a cycle $C_{\geq k+1}$ containing all vertices of $X$, a contradiction.

        \textbf{Case B.} $k$ is odd.

            Similarly, let $|X| = \frac{k+1}{2}$ for $X\subseteq D_+^G$, $Y = V(G)\setminus X$ and $G' = G(X, Y)$. Set $d'=d-\frac{k-1}{2}$. For each $x\in X$, $$\deg_{G'}(x) \geq \deg_G(x)-(|X|-1) \geq d-\frac{k-1}{2}=d'.$$ Since $3 \leq k < d$, we have $2 \leq |X|\leq d'$. Combined with $n \leq \frac{3}{2}d$, $$|Y| = n-|X| \leq \frac{3}{2}d-\frac{k+1}{2} = 2d-k-1-\frac{d-k-1}{2} \leq 2d-k-1 = 2d'-2.$$ Again Lemma~\ref{lemmaJackson1981} produces a cycle $C_{\geq k+1}$ in $G$ passing through all vertices of $X$, a contradiction.

        All parity cases yield a contradiction under $n \leq \tfrac{3}{2}d$. We conclude $n \geq \tfrac{3}{2}d + 1$, finishing the proof of Claim~\ref{claimLowerN}.
    \end{proof}\smallskip

    We next analyze the block decomposition of $G$. Assume first that the high-degree vertex set $D_+^G$ is not contained within a single 2-connected block of $G$. Then $G = G_1 \cup G_2$ with a single cut vertex $u = V(G_1)\cap V(G_2)$, so $n_1 + n_2 = n+1$ where $n_i = |V(G_i)|$. Denote $\theta_+^{G_i} = |\{v\in V(G_i)\mid \deg_{G_i}(v)\geq d\}|$. Since the cut vertex $u$ may lose high-degree status when restricted to each component, we have the upper bound
    $$
    \theta_+^G \leq \theta_+^{G_1} + \theta_+^{G_2} + 1.
    $$
    Define the threshold function
    $$
    \bar{\theta}(n_i) = \left\lfloor\frac{k}{2}\right\rfloor\cdot \frac{n_i-1}{d}
    +\frac{(n_i-d-1)\lfloor k/2\rfloor}{d(d-\lfloor k/2\rfloor)}.
    $$

    If either $\theta_+^{G_1}\geq \bar{\theta}(n_1)+1$ or $\theta_+^{G_2}\geq \bar{\theta}(n_2)+1$, then $G_1$ or $G_2$ forms a smaller counterexample violating vertex minimality Condition \eqref{psiItem3}. We thus assume $\theta_+^{G_i}\leq \bar{\theta}(n_i)$ for $i=1,2$. Substitute into the upper bound of $\theta_+^G$:
    \begin{align*}
        \theta_+^G &\leq \bar{\theta}(n_1) + \bar{\theta}(n_2) + 1 \\
        &= \left\lfloor\frac{k}{2}\right\rfloor\cdot \frac{n_1+n_2-2}{d}
        + \frac{(n_1+n_2-2d-2)\lfloor k/2\rfloor}{d(d-\lfloor k/2\rfloor)} + 1 \\
        &= \left\lfloor\frac{k}{2}\right\rfloor\cdot \frac{n-1}{d}
        + \frac{(n-2d-1)\lfloor k/2\rfloor}{d(d-\lfloor k/2\rfloor)} + 1 \\
        &< \left\lfloor\frac{k}{2}\right\rfloor\cdot \frac{n-1}{d}
        + \frac{(n-d-1)\lfloor k/2\rfloor}{d(d-\lfloor k/2\rfloor)} + 1 = \theta_+^G.
    \end{align*}
    The strict inequality $\theta_+^G < \theta_+^G$ is a direct contradiction. Therefore all vertices in $D_+^G$ must lie inside a unique 2-connected block of $G$.

    Let $G^*$ be a maximal $C_{\geq k+1}$-free supergraph in $G$, obtained by adding all possible edges with at least one endpoint in $D_+^G$ without creating a cycle $C_{\geq k+1}$. No edge between two vertices of $D_-^G$ is ever added. Observe that Claims \ref{claimNeighborOfLow} and \ref{claimLowerN} (structural properties concerning adjacency between vertices in $D_-^G$ and $D_+^G$) remain true in $G^*$, since we never add edges between two vertices of $D_-^G$. We split the rest of the proof into two exhaustive structural cases on the clique property of $D_+^G$ in $G^*$.

    \subsection{Case 1. $D_+^G$ forms a clique in $G^*$.}

        Let $C$ be a longest cycle in $G^*$ with the maximum number of high-degree vertices of $D_+^G$, with cyclic ordering $a_1,a_2,\dots,a_m$ along $C$. By the $C_{\geq k+1}$-free property, $m = |V(C)|\leq k$. Indices are taken modulo $m$ for all subscripts of vertices on $C$.

        \begin{claim2026ErdosQuestion}\label{claimDplusVc}
            Every high-degree vertex of $D_+^G$ lies on the longest cycle $C$, i.e., $D_+^G \subseteq V(C)$.
        \end{claim2026ErdosQuestion}
        \begin{proof}
            Suppose for contradiction there exists $u\in D_+^G \setminus V(C)$. We split two subcases based on whether adjacent vertices on $C$ are high-degree.
            \begin{itemize}
                \item  If every vertex of $C$ is a high-degree vertex of $D_+^G$, pick any consecutive pair $a_i,a_{i+1}\in V(C)$. The cycle $C' = a_i u a_{i+1} \overrightarrow{C} a_i$ is strictly longer than $C$, contradicting maximality of $C$.
                \item If some vertex $a_i\in V(C)$ is low-degree, Claim~\ref{claimNeighborOfLow}~\eqref{claimNeighborOfLow2} implies its two neighbors $a_{i-1},a_{i+1}\in D_+^G$. Construct a new cycle $C' = a_{i-1} u a_{i+1} \overrightarrow{C} a_{i-1}$. This cycle contains all vertices of $C$ except the single low-degree vertex $a_i$, hence $C'$ is also the longest cycle and carries more high-degree vertices than $C$, contradicting our selection rule for $C$.
            \end{itemize}
            Thus no such vertex $u$ exists, so $D_+^G\subseteq V(C)$.
        \end{proof}\smallskip

        \begin{claim2026ErdosQuestion}\label{claimNleq2dMinus2}
            The vertex set of $G^*$ satisfies $n \leq 2d - 2$.
        \end{claim2026ErdosQuestion}
        \begin{proof}
             By Claim~\ref{claimDplusVc}, all high-degree vertices of $D_+^G$ lie on $C$, so $\theta_+^G \leq |V(C)| \leq k$, otherwise, if $\theta_+^G \geq k+1$, then $C$ contains at least $k+1$ high-degree vertices, forcing a $C_{\geq k+1}$ cycle, a contradiction. Suppose for contradiction  $n\geq 2d-1$. We separate the proof by parity of $k$.

            \paragraph{Case A. $k$ is even.} \mbox{}\par

                Substitute $n = 2d - 1$ into the lower bound of $\theta_+^G$,
                \begin{align*}
                    \theta_+^G &\geq \frac{k}{2}\cdot \frac{2d - 1 - 1}{d}
                    + \frac{(2d-1-d-1)\lfloor k/2\rfloor}{d(d-\lfloor k/2\rfloor)} + 1 \\
                               &= k + 1 - \frac{k}{d} + \frac{k(d - 2)}{d(2d - k)} \\
                               &= k + 1 + \frac{k(k - d - 2)}{d(2d - k)}.
                \end{align*}
                We verify the fractional term satisfies $$\frac{k(k-d-2)}{d(2d-k)} > -1,$$ which rearranges to
                $$
                    k(k-d-2) > -d(2d-k) \iff k(k-2) + 2d(d-k) > 0.
                $$
                Since $2 \leq k < d$, both terms $k(k-2) \geq 0$ and $2d(d-k) > 0$, so the inequality holds. This yields $\theta_+^G > k$, contradicting $\theta_+^G \leq k$.

            \paragraph{Case B. $k$ is odd.} \mbox{}\par

                First assume $n \geq 2d$. Substitute $n=2d$ into the threshold,
                \begin{align*}
                    \theta_+^G &\geq \frac{k - 1}{2}\cdot \frac{2d - 1}{d}
                    + \frac{(2d-d-1)\lfloor k/2\rfloor}{d(d-\lfloor k/2\rfloor)} + 1 \\
                               &= k - \frac{k - 1}{2d} + \frac{(d - 1)(k - 1)}{d(2d - k + 1)} \\
                               &= k + \frac{(k - 3)(k - 1)}{2d(2d - k + 1)}.
                \end{align*}
                For odd $k > 3$, the additional fraction is positive, so $\theta_+^G > k$, contradicting $\theta_+^G \leq k$. If $k=3$, then $\theta_+^G \geq k.$ Since $\theta_+^G \leq k$, we get $\theta_+^G = k = 3$. By Claim~\ref{claimDplusVc}, $D_+^G \subseteq V(C)$, so $|V(C)| \geq 3$. But the $C_{\geq k+1}$-free property forces $|V(C)| \leq 3$. Hence $|V(C)| = \theta_+^G = 3$ and all vertices on $C$ are high-degree. So each vertex of $C$ has at least $d-2$ neighbors outside $C$. Any external vertex adjacent to two distinct vertices of $C$ creates a cycle $C_4 \subseteq C_{\geq k+1}$, a contradiction. Thus every vertex outside $C$ connects to at most one vertex on $C$. Since each vertex of $C$ has at least $d-2$ neighbors outside $C$, the number of vertices outside $C$ is at least $3(d-2)$. Count vertices $$n \geq 3(d - 2) + 3 = 3d - 3.$$ Substitute $n \geq 3d - 3$ into the threshold inequality,
                \begin{align*}
                    \theta_+^G \geq \frac{3d - 3 - 1}{d} + \frac{3d - 3 - d - 1}{d(d - 1)} + 1 = 4 - \frac{2}{d - 1}.
                \end{align*}
                Since $2< k < d$, we have $\frac{2}{d - 1} < 1$, and so $\theta_+^G > 3$. Since $\theta_+^G$ is an integer, $\theta_+^G \geq 4 > 3 = k$. By Claim~\ref{claimDplusVc}, there is a cycle $C_4$ in $G$, a contradiction.

                Next suppose $n = 2d - 1$. Substitute $n = 2d - 1$ into the threshold inequality,
                \begin{align*}
                    \theta_+^G &\geq \frac{k - 1}{2}\cdot \frac{2d - 1 - 1}{d}
                    + \frac{(2d-1-d-1)\lfloor k/2\rfloor}{d(d-\lfloor k/2\rfloor)} + 1 \\
                                &= k - \frac{k - 1}{d} + \frac{(d - 2)(k - 1)}{d(2d - k + 1)} \\
                                &= k + \frac{(k - 1)(k - d - 3)}{d(2d - k + 1)}.
                \end{align*}
                We now verify $$\frac{(k - 1)(k - d - 3)}{d(2d - k + 1)} > -1.$$ Since $2 \leq k < d$, the denominator $d(2d-k+1)>0$. Rearranging the inequality yields
                $$
                    (k-1)(k-d-3) + d(2d-k+1) > 0 \iff 2(d+1)(d-k) + k(k-2) + 3 > 0.
                $$
                Since $2 \leq k < d$, the left-hand-side is strictly positive. Hence
                $$
                    \theta_+^G \geq k+\frac{(k-1)(k-d-3)}{d(2d-k+1)} > k-1.
                $$
                Since $\theta_+^G$ is an integer, $\theta_+^G \geq k$. Combined with $\theta_+^G \leq k$, by Claim~\ref{claimDplusVc}, $\theta_+^G = |V(C)| = k$ and all vertices on $C$ are high-degree. For any vertex $u \in V(G^*)\setminus V(C)$, we have $|N(u)\cap V(C)|=1$. Otherwise, let $a_i, a_j \in N(u) \cap V(C)$, then $a_iua_j\overleftarrow{C}a_{i+1}a_{j+1}\overrightarrow{C}a_i$ is a cycle $C_{\geq k+1}$ in $G^*$ contains all vertices of $C$ and $u$, a contradiction. Each vertex on $C$ has at least $d-k+1$ neighbors outside $C$, hence
                $$
                    2d-1\geq k+k(d-k+1).
                $$
                Simplify the right-hand side minus left-hand side,
                $$
                    k+k(d-k+1)-(2d-1)=(d - k)(k - 2) + 1 > 0.
                $$
                This numerical contradiction shows $n \ne 2d - 1$.

                All possible values of $n$ lead to contradictions, hence $n \leq 2d - 2$.
        \end{proof}\smallskip

        \begin{claim2026ErdosQuestion}\label{claimVcK}
            The longest cycle $C$ satisfies $|V(C)| = k$.
        \end{claim2026ErdosQuestion}
        \begin{proof}
            By Claim~\ref{claimLowerN}~\eqref{claimItem1}, $\theta_+^G \geq \left\lfloor\frac{k}{2}\right\rfloor + 1$. Choose a subset $X\subseteq D_+^G$ with $|X|=\lceil{\frac{k+1}{2}}\rceil$, set $Y=V(G^*)\setminus X$, and let $G' = G^*(X,Y)$ be the spanning bipartite subgraph of $G^*$ with partitions $X$ and $Y$. Define $d'=d-\lceil{\frac{k-1}{2}}\rceil$. For each vertex $u$ of $X$, $$\deg_{G'}(u) \geq d-\left(\left\lceil\frac{k + 1}{2}\right\rceil-1\right) = d - \left\lceil\frac{k - 1}{2}\right\rceil = d'.$$

            Claim~\ref{claimNleq2dMinus2} gives $n\leq 2d - 2$. Note that $2 \leq k < d$. If $k$ is odd, then $|X| = \frac{k + 1}{2}$ and $d' = d - \frac{k-1}{2}$. Hence $2 \leq |X| \leq d'$ and
            \begin{align*}
                |Y| = n - |X|
                    \leq (2d - 2) - \frac{k + 1}{2}
                    \leq (2d - 2) + (d - k - 1) - \frac{k + 1}{2}
                    = 3d' - 5.
            \end{align*}
            If $k$ is even, then $|X|=\frac{k + 2}{2}$ and $d' = d - \frac{k}{2}$. Hence $2 \leq |X| \leq d'$. By Claim~\ref{claimDplusVc}, $D_+^G \subseteq V(C)$, so $|V(C)| \geq \theta_+^G$. If $n = 2d - 2$ and $k = d - 1$, then
            \begin{align*}
                |V(C)| \geq \theta_+^G
                       &\geq \left\lfloor\frac{k}{2}\right\rfloor\cdot \frac{n-1}{d}
                       + \frac{(n-d-1)\lfloor k/2\rfloor}{d(d-\lfloor k/2\rfloor)} + 1 \\
                       &= \frac{k(2k - 1)}{2(k + 1)} + 1 + \frac{k(k - 2)}{(k + 1)(k + 2)} \\
                       &= k - \frac{k - 2}{2k + 2} + \frac{k(k - 2)}{(k + 1)(k + 2)}.
            \end{align*}
            Since $\frac{k - 2}{2k + 2} < 1$ and $\frac{(k - 2)k}{(k + 1)(k + 2)} \geq 0$, we have $|V(C)| > k - 1$. Combine with $|V(C)| \leq k$, we get $|V(C)| = k$. The claim holds. If $n < 2d - 2$ or $k < d - 1$, then $$|Y| \leq n - |X| < (2d - 2) -  \frac{k + 2}{2} + (d - k - 1) \leq 3d - \frac{3}{2}k - 4 = 3d' - 4.$$ Since $|Y|$ is an integer, $|Y| \leq 3d' - 5.$ Thus in all cases $|Y| \leq 3d'-5$.

            If $G'$ is 2-connected, Lemma~\ref{lemma} gives a cycle $C_{\geq k+1}$ containing $X$, a contradiction.

            Now suppose $G'$ is 1-connected, and let $u$ be a cut vertex of $G'$. We first claim that $G$ contains at least two low-degree vertices. Suppose for contradiction that $G$ has only one low-degree vertex, say $v$. By Claim~\ref{claimNeighborOfLow}~\eqref{claimNeighborOfLow2}, every neighbor of $v$ is a high-degree vertex. By Claim~\ref{claimLowerN}~\eqref{claimItem1}, $\theta_+^G \geq \lfloor k/2 \rfloor + 1 \geq 2$, so there are at least two high-degree vertices, all of which are adjacent to $v$. Since $G$ has no isolated vertices, $v$ has at least one neighbor. Deleting $v$ leaves a graph in which every vertex has degree at least $d-1$, so by Dirac's theorem (Lemma~\ref{lemmaDirac}) it contains a cycle $C_{\geq d} \subseteq C_{\geq k+1}$, a contradiction. Hence $G$ contains at least two low-degree vertices. In particular, $G'$ contains a low-degree vertex $v \ne u$.

            Since $u$ is a cut vertex of $G'$, the graph $G'-u$ has at least two components. Let $C_1, C_2, \ldots, C_t$ ($t \geq 2$) be the components of $G'-u$. The low-degree vertex $v \ne u$ lies in some component, say $C_v$.

            If some component $C_j \ne C_v$ contains a vertex $a_i \in X$, then $C_v$ contains the low-degree vertex $v$ and $C_j$ contains $a_i \in X$, as desired.

            Otherwise, all vertices of $X$ lie in $C_v$. Then every other component $C_j$ ($j \ne v$) contains no vertex of $X$, so $C_j \subseteq Y$. Since $Y$ is an independent set in the bipartite graph $G'$, every vertex of $C_j$ is isolated in $G'-u$; hence $C_j$ consists of a single vertex $y_j \in Y$.

            If some $y_j$ is low-degree, then $C_j$ contains a low-degree vertex and $C_v$ contains some $a_i \in X$, as desired.

            It remains to consider the case where every $y_j$ is high-degree. Since all high-degree vertices of $G^*$ form a clique, each $y_j$ is adjacent in $G^*$ to every high-degree vertex in $C_v$. In particular, since $X \subseteq D_+^G$, each $y_j$ is adjacent to every vertex of $X$. But $X$ is nonempty and $X \subseteq C_v$, so $y_j$ is adjacent to some vertex of $X$ in $G^*$. This edge lies between $X$ and $Y$, hence it belongs to $G'$, contradicting the fact that $y_j$ is isolated in $G'$. Therefore this case cannot occur.

            Thus there exists a component containing a low-degree vertex $v \in D_-^G$ and another component containing some vertex $a_i \in X$.

            Suppose for contradiction that every $v$-$a_i$ path in $G^*$ has at most $k$ vertices. The edge $va_i$ has one endpoint in $D_+^G$, hence it is eligible under the maximality rule of $G^*$. Adding $va_i$ can only create cycles formed by concatenating $va_i$ with some $v$-$a_i$ path in $G^*$, all of which have at most $k$ vertices. Thus adding $va_i$ cannot produce any $C_{\geq k+1}$ cycles, contradicting the maximality of $G^*$.

            Therefore there exists a $v$-$a_i$ path $P$ in $G^*$ with $|V(P)|\geq k+1$. Write $P=vb_1b_2\cdots b_\ell a_i$, so that $\ell+2=|V(P)|\geq k+1$, i.e., $\ell+1\geq k$. By Claim~\ref{claimNeighborOfLow}~\eqref{claimNeighborOfLow2}, every neighbor of the low-degree vertex $v$ is a high-degree vertex; in particular, $b_1\in D_+^G$. Since $a_i\in X\subseteq D_+^G$ and all high-degree vertices of $G^*$ form a clique, the vertices $b_1$ and $a_i$ are adjacent in $G^*$. Hence $b_1b_2\cdots b_\ell a_i b_1$ is a cycle in $G^*$ of length $\ell+1\geq k$. Since $C$ is a longest cycle of $G^*$ and $|V(C)|\leq k$, we obtain $k\leq |V(C')|\leq |V(C)|\leq k$, and therefore $|V(C)|=k$.
        \end{proof}\smallskip

        Let
        \begin{align*}
            X=\{a_i \in D_+^G | a_{i-1}\not\in D_+^G\text{ and }a_{i+1}\not\in D_+^G\},\qquad
            Y_1=\{a_i \in D_+^G | a_{i-1}\in D_+^G\text{ and }a_{i+1}\not\in D_+^G\},\\
            Y_2=\{a_i \in D_+^G | a_{i-1}\not\in D_+^G\text{ and }a_{i+1}\in D_+^G\},\qquad
            Z=\{a_i \in D_+^G | a_{i-1}\in D_+^G\text{ and }a_{i+1} \in D_+^G\},
        \end{align*}
        and write $$|X| = x,\ \ |Y_1| = |Y_2| = y,\ \ |Z| = z.$$ We have
        \begin{equation}\label{equationX2yd}
            \begin{split}
                x+2y+z = \theta_+^G,
            \end{split}
        \end{equation}
        Count boundary transitions between high/low vertices on $C$,
        $$
            |\{a_i\in D_+ \mid a_{i+1}\notin D_+^G\}|=x+y.
        $$
        By Claim~\ref{claimVcK}, $|V(C)|=k$, the number of low segments cannot exceed total remaining vertices,
        \begin{equation}\label{equationXy}
            \begin{split}
                 x+y \leq k - \theta_+^G.
            \end{split}
        \end{equation}
        Combine \eqref{equationX2yd} with \eqref{equationXy},
        \begin{equation}\label{equationYz}
            \begin{split}
               y+z \geq  2\theta_+^G - k.
            \end{split}
        \end{equation}

        Let $S$ denote the collection of paths along $C$ connecting vertices in $Y_1$ to vertices in $Y_2$.

        \begin{claim2026ErdosQuestion}\label{claimY1Y2Z}

            Take two distinct vertices $a_i,a_j\in Y_1\cup Y_2\cup Z$ sharing a common external neighbor $u\notin V(C)$. Then
            \begin{enumerate}[(i)]
                \item $a_i,a_j$ are exactly the two endpoints of some path in $S$;\label{Y1Y2ZItem1}
                \item If $a_i\in Y_1$ and $a_j\in Y_2$, then neither $a_{i+1}$ nor $a_{j-1}$ has edges to any vertex in $(Y_1\cup Y_2\cup Z)\setminus\{a_i,a_j\}$.\label{Y1Y2ZItem2}
            \end{enumerate}
        \end{claim2026ErdosQuestion}
        \begin{proof}
            Let $u\notin V(C)$ satisfy $ua_i,ua_j\in E(G^*)$.

            \emph{(i)} Assume $a_i,a_j$ are not endpoints of some path in $S$.
            \begin{itemize}
                \item If $a_i,a_j\in Y_1\cup Z$, then $a_{i-1},a_{j-1}\in D_+^G$. Hence $a_{i-1}a_{j-1}\in E(G^*)$ and $C'=a_i \overrightarrow{C} a_{j-1} a_{i-1} \overleftarrow{C} a_j u a_i$ is a cycle in $G^*$ containing $u$ and all vertices of $C$, which is longer than $C$, a contradiction.
                \item If $a_i,a_j\in Y_2\cup Z$, then $a_{i+1},a_{j+1}\in D_+^G$. Hence $a_{i+1}a_{j+1}\in E(G^*)$ and $C'=a_i \overleftarrow{C} a_{j+1} a_{i+1} \overrightarrow{C} a_j u a_i$ is a cycle in $G^*$ containing $u$ and all vertices of $C$, which is longer than $C$, a contradiction.
                 \item If $a_i\in Y_1,a_j\in Y_2$ and not paired in $S$, the segment $a_{i+1}\overrightarrow{C}a_{j-1}$ contains some $a_t\in Y_1 \cup Y_2$; otherwise, $a_i,a_j$ are the two endpoints of some path in $S$ and \eqref{Y1Y2ZItem1} holds. Without loss of generality, assume $a_t\in Y_2$. Then $a_t,a_{t+1}\in D_+^G$. Hence $a_{i-1}a_t,a_{j+1}a_{t+1}\in E(G^*)$ and $C'=a_i\overrightarrow{C}a_ta_{i-1}\overleftarrow{C}a_{j+1}a_{t+1}\overrightarrow{C}a_jua_i$ is a cycle in $G^*$ containing $u$ and all vertices of $C$, which is longer than $C$, a contradiction.
            \end{itemize}
            By symmetry the case $a_i\in Y_2,a_j\in Y_1$ also fails. So \eqref{Y1Y2ZItem1} holds.

            \emph{(ii)} Let $a_i\in Y_1$ and $a_j\in Y_2$. By \eqref{Y1Y2ZItem1}, $a_i$ and $a_j$ are the two endpoints of some path in $S$. So $\{a_{i-1},a_i,a_j,a_{j+1}\}\subseteq D_+^G$. Suppose there exists $a_r\in (Y_1\cup Y_2\cup Z)\setminus\{a_i,a_j\}$ such that $a_{j-1}a_r\in E(G^*)$.
            \begin{itemize}
                \item If $a_r\in Y_1\cup Z$, then $a_{r-1}\in D_+^G$. Hence $a_{i-1}a_{r-1}\in E(G)$ and $C'=a_i\overrightarrow{C}a_{j-1}a_r\overrightarrow{C}a_{i-1}a_{r-1}\overleftarrow{C}a_jua_i$ is a cycle in $G^*$ containing $u$ and all vertices of $C$, which is longer than $C$, a contradiction.
                \item If $a_r\in Y_2$, then $a_{r+1}\in D_+^G$. Hence $a_{j+1}a_{r+1}\in E(G^*)$ and $C'=a_i\overrightarrow{C}a_{j-1}a_r\overleftarrow{C}a_{j+1}a_{r+1}\overrightarrow{C}a_{i-1}a_jua_i$ is a cycle in $G^*$ containing $u$ and all vertices of $C$, which is longer than $C$, a contradiction.
            \end{itemize}
            Identical argument applies to $a_{i+1}$, completing \emph{(ii)}.
        \end{proof}\smallskip

        Let $t=|Y_1'|=|Y_2'|$, where $Y_1'\subseteq Y_1,Y_2'\subseteq Y_2$ are the vertices paired with shared external vertices, $t\leq y$. By Claims \ref{claimVcK} and Claim~\ref{claimY1Y2Z}~\eqref{Y1Y2ZItem2},
        \begin{itemize}
            \item Vertices in $(Y_1\cup Y_2\cup Z)\setminus(Y_1'\cup Y_2')$ connect to at most $k-1-t$ vertices of $C$;
            \item Vertices in $Y_1'\cup Y_2'$ connect to at most $k-t$ vertices of $C$.
        \end{itemize}
        Count total vertices of $G^*$,
        $$
            \begin{aligned}
                n &\geq |V(C)| + |(Y_1 \cup Y_2 \cup Z) \backslash \{Y_1',Y_2'\}|(d-(k-1-t)) + |Y_1'|(d-(k-t)) \\
                    &=k+(2y+z-2t)(d-k+1+t)+t(d-k+1+t)-t,
            \end{aligned}
        $$
        rearranged to
        \begin{equation}\label{equationNMT1}
            \begin{split}
               n-k+t \geq (2y+z-t)(d-k+1+t).
            \end{split}
        \end{equation}
        Since $t\leq y$, combine \eqref{equationYz} with \eqref{equationNMT1},
        \begin{equation}\label{equationNMT2}
            \begin{split}
                n-k+t \geq (2\theta_+^G-k)(d-k+1+t)
                      = (2\theta_+^G -k)(d-k+1) + (2\theta_+^G -k-1)t.
            \end{split}
        \end{equation}
        Drop the nonnegative term $(2\theta_+^G -k-1)t$,
        \begin{equation}\label{equationNMinusK}
            \begin{split}
                n-k \geq (2\theta_+^G -k)(d-k+1).
            \end{split}
        \end{equation}
        By Claim~\ref{claimLowerN}~\eqref{claimItem1}, substitute the lower bound of $\theta_+^G$ into \eqref{equationNMinusK},
        \begin{equation}\label{equationNM1}
            \begin{split}
                n-k \geq \left(2\left(\left\lfloor{\frac{k}{2}}\right\rfloor \cdot \frac{n-1}{d}
                +\frac{(n-d-1)\lfloor k/2\rfloor}{d(d-\lfloor k/2\rfloor)}+1\right)-k\right)(d-k+1).
            \end{split}
        \end{equation}
        Recall Claim~\ref{claimNleq2dMinus2}: $n\leq 2d-2$, so
        \begin{equation}\label{equation2d2knk}
            \begin{split}
               2d-2-k\geq n-k.
            \end{split}
        \end{equation}

        Recall Claim~\ref{claimLowerN}~\eqref{claimItem2}: $n\geq \tfrac32 d+1$. Combine \eqref{equationNM1} with \eqref{equation2d2knk} based on the parity of $k$.

        If $k$ is even, drop the positive term $\frac{2k(n - d - 1)}{d(2d - k)}$, then
        $$
               2d - 2 - k \geq (\frac{k(n - 1)}{d} + \frac{2k(n - d - 1)}{d(2d - k)} + 2 - k) (d - k + 1) \\
                      > (\frac{k}{2} + 2) (d - k + 1).
        $$
        But direct expansion, since $d > k$, $$(\frac{k}{2} + 2) (d - k + 1) - (2d - 2 - k) = \frac{k(d - k - 1) + 8}{2} > 0,$$ a numerical contradiction.

        If $k$ is odd, drop the positive term $\frac{2(k - 1)(n - d - 1)}{d(2d - k + 1)} > 0$, then
        \begin{equation}\label{equationNK2Odd}
            \begin{split}
               2d - 2 - k &\geq (\frac{(k - 1)(n - 1)}{d} + \frac{2(k - 1)(n - d - 1)}{d(2d - k + 1)} + 2 - k) (d - k + 1) \\
                      &> \frac{k + 1}{2} (d - k + 1).
            \end{split}
        \end{equation}
        Set $\beta=\tfrac{k+1}{2}(d-k+1)-(2d-2-k)$, write $d=k+x$. By \eqref{equationNK2Odd}, $\beta < 0$.
        \begin{equation}\label{equationNK2Odd2}
            \begin{split}
               \beta = \frac{(x - 1)k - 3x + 5}{2}.
            \end{split}
        \end{equation}
       Since $d > k$, we have $x \geq 1$. We analyze all positive integer possibilities for $x\geq 1$.
        \begin{itemize}
            \item If $x=1$:
                $$
                \beta=\frac{(1-1)k-3\cdot 1+5}{2}=1>0.
                $$
            \item If $x=2$:
                $$
                \beta=\frac{(2-1)k-3\cdot 2+5}{2}=\frac{k-1}{2}.
                $$
                Since $k\geq 3$, we have $\frac{k-1}{2}\geq 1$, hence $\beta>0$.
            \item If $x=3$:
                $$
                \beta=\frac{(3-1)k-3\cdot 3+5}{2}=\frac{2k-4}{2}=k-2.
                $$
                Since $k\geq 3$, we get $\beta = k-2 \geq 1>0$.
            \item If $x\geq 4$: using $k\geq 3$,
                $$
                (x-1)k-3x+5 \geq 3(x-1)-3x+5 = 2>0,
                $$
                which implies $\beta>0$.
        \end{itemize}

        Thus for every integer $x\geq 1$, we obtain $\beta>0$. This stands in direct contradiction to $\beta<0$ coming from inequality \eqref{equationNK2Odd}. All subcases yield contradictions.

\subsection{Case 2. $D_+^G$ is not a clique in $G^*$.}

In this case, there exist two vertices in $D_+^G$ that are not joined by an edge of $G^*$. Since adding an edge between two such vertices would create a $C_{\geq k+1}$ cycle, they must be connected by a path of length at least $k+1$. Let $P=a_1a_2\cdots a_m$ be a longest path in $G^*$ with $a_1,a_m\in D_+^G$. Then $m\geq k+1.$ Give $P$ an orientation from $a_1$ to $a_m$ denoted by $\overrightarrow{P}$, and the reverse orientation denoted by $\overleftarrow{P}$. For $w\in\{1,m\}$, define
$$
    Q_w=N_{G^*}(a_w)\setminus V(P),\qquad
    R_w=N_{G^*}(Q_w)\cap D_+^G\cap V(P),\qquad
    S_w=N_{G^*}(a_w)\cap V(P)\setminus R_w,
$$
and let
$$
    q_w=|Q_w|,\qquad r_w=|R_w|,\qquad s_w=|S_w|.
$$
For any $U\subseteq V(P)$, let
$$
    U^-=\{a_j\in V(P)\mid a_{j+1}\in U\},\qquad
    U^+=\{a_j\in V(P)\mid a_{j-1}\in U\}.
$$
Set
$$
    t_1=|S_1\setminus R_1^-|,\qquad t_m=|S_m^+\setminus R_m^{++}|.
$$

If $Q_w=\emptyset$, then all neighbors of $a_w$ lie on $P$. Since $\deg_{G^*}(a_w)\geq d>k$, pick the farthest neighbor $a_i\in N(a_w)$; then $a_w\overrightarrow{P}a_i a_w$ is a $C_{\geq k+1}$ cycle, a contradiction. Thus $Q_w\neq\emptyset$. Also $a_w\in R_w$, so
$$
    q_w+r_w+s_w\geq |N_{G^*}(a_w)|+|\{a_w\}|\geq d+1,
$$
giving
\begin{equation}\label{eq:rw-lower}
    r_w\geq d+1-s_w-q_w.
\end{equation}
Maximality of $P$ implies
\begin{equation}\label{eq:Qw-subset}
    Q_w\subseteq D_-^G\quad\text{and}\quad N_{G^*}(Q_w)\cap D_+^G=R_w.
\end{equation}
Assume
$$
    r_w\leq \left\lfloor\frac{k}{2}\right\rfloor\cdot\frac{q_w}{d}+\frac{q_w\lfloor k/2\rfloor}{d(d-\lfloor k/2\rfloor)}.
$$
Since $Q_w\subseteq D_-^G$, deleting the vertices of $Q_w$ does not remove any vertex from $D_+^G$. Hence the number of original high-degree vertices in $G-Q_w$ is still $\theta_+^G$. Since
\begin{align*}
    &\left\lfloor\frac{k}{2}\right\rfloor\cdot\frac{n-q_w-1}{d}+\frac{(n-q_w-d-1)\lfloor k/2\rfloor}{d(d-\lfloor k/2\rfloor)}+1\\
    &=\left\lfloor\frac{k}{2}\right\rfloor\cdot\frac{n-1}{d}+\frac{(n-d-1)\lfloor k/2\rfloor}{d(d-\lfloor k/2\rfloor)}+1
    -\left(\left\lfloor\frac{k}{2}\right\rfloor\cdot\frac{q_w}{d}+\frac{q_w\lfloor k/2\rfloor}{d(d-\lfloor k/2\rfloor)}\right)\\
    &\leq \theta_+^G-r_w,
\end{align*}
the graph $G-Q_w$ has at least
$$
    \left\lfloor\frac{k}{2}\right\rfloor\cdot\frac{n-q_w-1}{d}+\frac{(n-q_w-d-1)\lfloor k/2\rfloor}{d(d-\lfloor k/2\rfloor)}+1
$$
high-degree vertices and is $C_{\geq k+1}$-free with fewer vertices, contradicting the minimality of $G$. Therefore
$$
    r_w>\left\lfloor\frac{k}{2}\right\rfloor\cdot\frac{q_w}{d}+\frac{q_w\lfloor k/2\rfloor}{d(d-\lfloor k/2\rfloor)},
$$
that is,
$$
    \frac{dr_w}{\lfloor k/2\rfloor\left(1+\frac{1}{d-\lfloor k/2\rfloor}\right)}>q_w.
$$
Define
$$
    \alpha(d)=\frac{d}{\lfloor k/2\rfloor\left(1+\frac{1}{d-\lfloor k/2\rfloor}\right)}=\frac{d(d-\lfloor k/2\rfloor)}{\lfloor k/2\rfloor(d-\lfloor k/2\rfloor+1)}.
$$
Thus
\begin{equation}\label{eq:alpha-rw}
    \alpha(d)r_w>q_w.
\end{equation}

By maximality of $P$, $R_w\cap R_w^-=\emptyset$ and $R_m^+\cap R_m^{++}=\emptyset$. By definition $S_w\cap R_w=\emptyset$. Since $m\geq k+1$, any non-empty intersection among $S_1,R_1,R_1^-,S_m^+,R_m^+,R_m^{++}$ other than $S_1\cap R_1^-$ and $S_m^+\cap R_m^{++}$ forces a $C_{\geq k+1}$ in $G^*$, a contradiction. Since $a_1\in R_1$, $|R_1^-|=|R_1\setminus\{a_1\}|=r_1-1$. Maximality gives $a_{m-1}\notin R_m$, hence $|R_m^{++}|=|R_m^+|=|R_m\setminus\{a_m\}|=r_m-1$. From set-size estimates
\begin{equation}\label{eq:rw-sw}
    r_w\geq s_w+1-t_w.
\end{equation}
Combine \eqref{eq:rw-lower} with \eqref{eq:alpha-rw} and \eqref{eq:rw-sw} to eliminate $s_w,q_w$:
\begin{equation}\label{eq:rw-B}
    r_w>\frac{d+2-t_w}{\alpha(d)+2}.
\end{equation}
Define
$$
    X=R_1\cup R_1^-\cup S_1\cup R_m^+\cup R_m^{++}\cup S_m^+.
$$

\begin{claim2026ErdosQuestion}\label{claim:Xcycle}
    There exists a cycle $C_{|X|-1}$ in $G^*$.
\end{claim2026ErdosQuestion}

\begin{proof}
The construction is illustrated in Figures \ref{figureClaim7Even} and \ref{figureClaim7Odd}. We construct a sequence of disjoint paths $P_1, P_2, \ldots, P_h$ iteratively as follows.

Let $P_1$ be a path outside $P$ (except for its end vertices) connecting $a_{x_1}(=a_1)$ to $a_{y_1}$ such that $y_1$ is as large as possible. Assume that we have constructed $P_1, P_2, \ldots, P_f$, where $P_f$ connects $a_{x_f}$ and $a_{y_f}$. If $a_{y_f} = a_m$, then we have done; otherwise, let $P_{f+1}$ be a path outside $P$ (except for its end vertices) connecting $a_{x_{f+1}}$ to $a_{y_{f+1}}$ such that $y_{f+1}$ is as large as possible subject to $x_{f+1} < y_f$.

By the 2-connectedness of $B$, we must have $y_{f+1} > y_f$. Otherwise, since $x_{f+1} < y_f$, the path $P_{f+1}$ connects two vertices on the segment $a_1 \overrightarrow{P} a_{y_f}$ (possibly including $a_{y_f}$ itself). This would imply that any path in $B$ from $a_1 \overrightarrow{P} a_{y_f}$ to $a_{y_f} \overrightarrow{P} a_m$ must pass through $a_{y_f}$, making $a_{y_f}$ a cut vertex of $B$, a contradiction.

Furthermore, we must have $x_{f+1} \geq y_{f-1}$. Otherwise, if $x_{f+1} < y_{f-1}$, then since $y_{f+1} > y_f$, the path $P_{f+1}$ provides a valid path from $a_{x_{f+1}}$ to $a_{y_{f+1}}$ with $x_{f+1} < y_{f-1}$ and $y_{f+1} > y_f$. This contradicts the choice of $P_f$, which requires $y_f$ to be as large as possible among all paths connecting a vertex before $y_{f-1}$ to a vertex after $y_f$.

Moreover, $P_{f+1}$ is disjoint from all paths $P_1, P_2, \ldots, P_f$ except for the possibility that $x_{f+1} = y_{f-1}$. Otherwise, if there is an intersection outside $P$ between $P_{f+1}$ and some path $P_g$ ($1 \leq g \leq f$), then we can replace $y_g$ with $y_{f+1}$ and re-construct $P_g, P_{g+1}, \ldots$, contradicting the maximality of the $y$-indices.

Since the indices $y_f$ strictly increase at each step and $P$ is finite, the construction must terminate after finitely many steps. It can only terminate when $a_{y_h} = a_m$, as otherwise we could construct another path $P_{h+1}$ with $y_{h+1} > y_h$.

Change $y_1$ and $x_h$ if necessary such that $y_1$ is as small as possible and $x_h$ is as large as possible without violating $x_2 < y_1$ and $x_h < y_{h-1}$, respectively. So, no vertex of $R_1$ lies in $a_{x_2+1}\overrightarrow{P}a_{y_1-1}$, and no vertex of $R_m \cup S_m$ lies in $a_{x_h+1}\overrightarrow{P}a_{y_{h-1}-1}$.

\begin{figure}
   \centering\includegraphics[width=9cm,height=2cm]{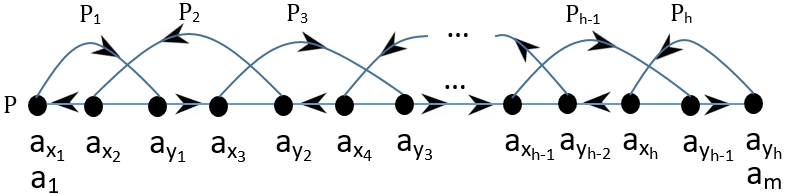}
   \caption{A cycle constructed for the case when $h$ is even.}
   \label{figureClaim7Even}
\end{figure}
\begin{figure}
   \centering\includegraphics[width=10cm,height=2cm]{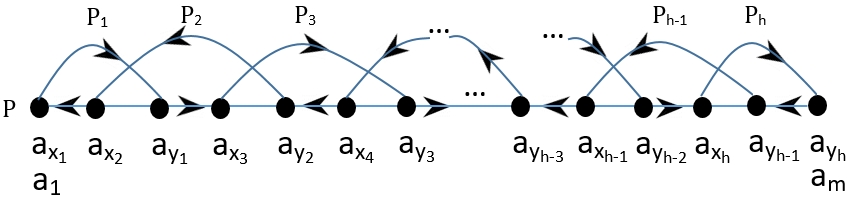}
   \caption{A cycle constructed for the case when $h$ is odd.}
   \label{figureClaim7Odd}
\end{figure}

As shown in Figures \ref{figureClaim7Even} and \ref{figureClaim7Odd}, if $h$ is even, then $G^*$ has a cycle $$C=(a_1=)a_{x_1}\overrightarrow{P_1}a_{y_1}\overrightarrow{P}a_{x_3}\overrightarrow{P_3}a_{y_3}\overrightarrow{P}\ldots a_{x_{h-1}}\overrightarrow{P_{h-1}}a_{y_{h-1}}\overrightarrow{P}a_{y_h}(=a_m)\overleftarrow{P_h}a_{x_h}\overleftarrow{P}a_{y_{h-2}}\ldots a_{x_4}\overleftarrow{P}a_{y_2}\overleftarrow{P_2}a_{x_2}\overleftarrow{P}a_{x_1};$$ and if $h$ is odd, then $G^*$ has a cycle $$C=(a_1=)a_{x_1}\overrightarrow{P_1}a_{y_1}\overrightarrow{P}a_{x_3}\overrightarrow{P_3}a_{y_3}\overrightarrow{P}\ldots a_{y_{h-2}}\overrightarrow{P}a_{x_h}\overrightarrow{P_h}a_{y_h}(=a_m)\overleftarrow{P}a_{y_{h-1}}\overleftarrow{P_{h-1}}a_{x_{h-1}}\overleftarrow{P}a_{y_{h-3}}\ldots a_{x_4}\overleftarrow{P}a_{y_2}\overleftarrow{P_2}a_{x_2}\overleftarrow{P}a_{x_1}.$$ In both cases, $C$ contains all vertices of $P$ except $a_{x_z+1}\overrightarrow{P}a_{y_{z-1}-1}$, where $z=2,3,\ldots,h$. Note that no vertex of $a_{x_3}\overrightarrow{P}a_m$, except $a_{y_1}$ when $x_3 = y_1$, is in $R_1 \cup S_1$ since $x_3 \geq y_1$; and no vertex of $a_1\overrightarrow{P}a_{y_{h-2}}$, except $a_{x_h}$ when $y_{h-2} = x_h$, is in $R_m \cup S_m$ since $y_{h-2} \leq x_h$. So, $C$ contains all vertices of $X$ except at most $a_{y_1-1}$, $a_{x_h+1}$ and $a_{x_h+2}$, which could be in $R_1^-$, $R_m^+ \cup S_m^+$ and $R_m^{++}$, respectively. But if $a_{y_1-1} \in R_1^-$, then $P_1$ must contain at least a vertex outside $P$. Similarly, if $a_{x_h+1} \in R_m^+$, then $P_h$ must contain at least a vertex outside $P$. So, $C$ is a cycle $C_{|X|-1}$ in $G^*$.
\end{proof}

Let $A_1=R_1\cup R_1^-\cup S_1$ and $A_2=R_m^+\cup R_m^{++}\cup S_m^+$. The sets $A_1,A_2$ are disjoint, and
$$
    |A_1|=2r_1-1+t_1,\qquad |A_2|=2r_m-2+t_m,\qquad |X|=2r_1+2r_m-3+t_1+t_m.
$$
By \eqref{eq:rw-B} and Claim~\ref{claim:Xcycle},
\begin{align}\label{eq:main-ineq-general}
    k+1 & \geq 2r_1+2r_m-3+t_1+t_m \nonumber\\
        & > 2\cdot\frac{d+2-t_1}{\alpha(d)+2}+2\cdot\frac{d+2-t_m}{\alpha(d)+2}-3+t_1+t_m.
\end{align}

\begin{definition2026ErdosQuestion}[$Q_1'$-technique]\label{def:Q1prime}
Assume $t_1=0$, i.e., $S_1\subseteq R_1^-$. Then $a_2\notin R_1$. Indeed, if $a_2\in R_1$, then there exists $u\in Q_1$ such that $ua_2\in E(G^*)$. Then $a_1ua_2\overrightarrow{P}a_m$ is an $a_1$-$a_m$ path longer than $P$, contradicting the maximality of $P$. Hence $a_2\notin R_1$. Since $a_2\in N_{G^*}(a_1)\cap V(P)$, we have $a_2\in S_1$. Because $S_1\subseteq R_1^-$, it follows that $a_2\in R_1^-$, and therefore $a_3\in R_1$. So there exists $a_2'\in Q_1$ with $a_2'a_3\in E(G^*)$. We claim that $a_2$ is low-degree, i.e., $a_2\in D_-^G$. Otherwise, if $a_2\in D_+^G$, then $a_2a_1a_2'a_3\overrightarrow{P}a_m$ is a path with endpoints in $D_+^G$ (namely $a_2$ and $a_m$) containing $a_2'$ and all vertices of $P$, which has $m+1$ vertices and is therefore longer than $P$, contradicting the maximality of $P$. Hence $a_2\in D_-^G$.

Define
$$
    Q_1'=Q_1\cup\{a_2\},\qquad
    R_1'=N_{G^*}(Q_1')\cap D_+\cap V(P),\qquad
    S_1'=N_{G^*}(a_1)\cap V(P)\setminus R_1',
$$
with
$$
    q_1'=|Q_1'|,\qquad r_1'=|R_1'|,\qquad s_1'=|S_1'|,\qquad
    t_1'=|S_1'\setminus R_1'^-|.
$$
Set
$$
    X'=R_1'\cup R_1'^-\cup S_1'\cup R_m^+\cup R_m^{++}\cup S_m^+.
$$
Since $Q_1'=Q_1\cup\{a_2\}$, the combinatorial construction argument in Claim~\ref{claim:Xcycle} applies to the triple $(Q_1',R_1',S_1')$ associated with $P$, provided that no vertex of $a_{x_3}\overrightarrow{P}a_m$ belongs to $R_1'\cup S_1'$. If such an interfering vertex $\bar a$ exists, we interchange $a_2'$ and $a_2$, replacing $P$ by the longest path $P'=a_1a_2'a_3\overrightarrow{P}a_m$. Iterating this interchange procedure, we can always choose a longest-path representative so that no such interfering vertex $\bar a$ occurs. For this representative, the construction of Claim~\ref{claim:Xcycle} yields a cycle $C_{|X'|-1}$ in $G^*$. Since $G^*$ is $C_{\geq k+1}$-free, we have
\begin{equation}\label{eq:Xprimeupper}
    |X'|\leq k+1.
\end{equation}
Since $|X'|=2r_1'+2r_m-3+t_1'+t_m$, we obtain
\begin{equation}\label{eq:r1primermupper}
    r_1'+r_m\leq \frac{k+4-(t_1'+t_m)}{2}\leq \frac{k+4}{2}.
\end{equation}
On the other hand, combine \eqref{eq:rw-lower} with \eqref{eq:rw-sw} to eliminate $s_w$:
$$
q_1\geq d+2-2r_1-t_1=d+2-2r_1.
$$
Since $q_1'=q_1+1$, $q_1'\geq d+3-2r_1$. By \eqref{eq:alpha-rw} applied to $Q_1'$, $\alpha(d)r_1'>q_1'$, hence
\begin{equation}\label{eq:r1primelower}
    r_1'>\frac{d+3-2r_1}{\alpha(d)}.
\end{equation}
We call this construction and the resulting inequalities \eqref{eq:Xprimeupper}--\eqref{eq:r1primelower} the \emph{$Q_1'$-technique}.
\end{definition2026ErdosQuestion}

Let $T=t_1+t_m\geq 2$.

\subsubsection{Case 2.1. $T=0$.}

In this case, $S_1\subseteq R_1^-$ and $S_m^+\subseteq R_m^{++}$. Since $T=0$, $t_1=t_m=0$. By Claim~\ref{claim:Xcycle},
\begin{equation}\label{eq:r1rm-upper}
    2r_1+2r_m-3+t_1+t_m\leq k+1\implies r_1+r_m\leq\frac{k+4}{2}.
\end{equation}

Define
$$
    B(d)=\frac{d+2}{\alpha(d)+2}.
$$
With $t_1=t_m=0$, \eqref{eq:rw-B} becomes $r_w>B(d)$. Since $r_w$ is an integer,
\begin{equation}\label{eq:rw-integer}
    r_w\geq \lfloor B(d)\rfloor+1.
\end{equation}
Together with \eqref{eq:r1rm-upper},
\begin{equation}\label{eq:b-upper}
    \lfloor B(d)\rfloor+1\leq\frac{k+4}{4}.
\end{equation}
If $\lfloor B(d)\rfloor+1>\frac{k+4}{4}$, we obtain a contradiction.

\begin{claim2026ErdosQuestion}\label{claim:Bmonotone}
    For $k\geq 6$, $B(d)$ is strictly increasing over integers $d>k$.
\end{claim2026ErdosQuestion}

\begin{proof}
Define $k'=\lfloor k/2\rfloor$ and the finite difference
$$
    \Delta B(d)=B(d+1)-B(d)=\frac{d+3}{\alpha(d+1)+2}-\frac{d+2}{\alpha(d)+2}.
$$
The denominator $(\alpha(d)+2)(\alpha(d+1)+2)$ is strictly positive since $\alpha(d)>0$. Hence $\Delta B(d)>0$ is equivalent to positivity of
$$
    \mathrm{Num}(d)=(d+3)\alpha(d)-(d+2)\alpha(d+1)+2.
$$
Substitute
$$
    \alpha(d)=\frac{d(d-k')}{k'(d-k'+1)},\qquad
    \alpha(d+1)=\frac{(d+1)(d+1-k')}{k'(d+2-k')}.
$$
Putting everything over the common denominator $k'(d-k'+1)(d+2-k')>0$, positivity of $\mathrm{Num}(d)$ is equivalent to positivity of
$$
    \mathcal{N}(d)=(d+3)d(d-k')(d+2-k')-(d+2)(d+1)(d+1-k')(d-k'+1)
$$
$$
    +2k'(d-k'+1)(d+2-k').
$$
Set $z=d-k'$. Expanding and collecting terms,
$$
    \mathcal{N}(z,k')=(2k'-3)z^2+(4k'-7)z-(k'^2-k'+2).
$$
For $k\geq 6$, $k'\geq 3$. The quadratic has positive leading coefficient $2k'-3>0$ and negative constant term $-(k'^2-k'+2)$. By Vieta's formulas, the product of its roots is negative, so it has exactly one positive real root $z_0(k')$ and one negative root. The vertex is at negative $z$, so $\mathcal{N}(z,k')$ is strictly increasing for $z>0$. Since $\mathcal{N}(k',k')=2k'^3-6k'-2>0$ for $k'\geq 3$, we have $z_0(k')<k'$. For admissible parameters, if $k$ is even then $d>2k'$ and $z=d-k'>k'$; if $k$ is odd then $d>2k'+1$ and $z>k'+1$. In both cases $z>k'>z_0(k')$, so $\mathcal{N}(z,k')>0$, hence $\Delta B(d)>0$. Thus $B(d)$ is strictly increasing over integers $d>k$, $k\geq 6$.
\end{proof}

\paragraph{Case 2.1.1. $k$ is even.}\mbox{}\par

Then
$$
    \alpha(d)=\frac{2d(2d-k)}{k(2d-k+2)}.
$$
Substituting the definition of $\alpha(d)$ into \eqref{eq:b-upper}, we analyze the following subcases.

First assume $k\equiv 2\pmod 4$.

Write $k=4x+2$ with integer $x\geq 1$. For $d=k+1=4x+3$, we have
$$
    \alpha(k+1)=\frac{2(4x+3)(x+1)}{(2x+1)(2x+3)}, \qquad \alpha(k+1)+2=\frac{16x^2+30x+12}{(2x+1)(2x+3)}.
$$
Therefore
$$
    B(k+1)=\frac{k+3}{\alpha(k+1)+2}=\frac{(4x+5)(2x+1)(2x+3)}{16x^2+30x+12}.
$$
We verify that $B(k+1)>x+1$. Since the denominator $16x^2+30x+12$ is positive, we have
\begin{align*}
    (4x+5)(2x+1)(2x+3)-(x+1)(16x^2+30x+12)&=(16x^3+52x^2+52x+15)-(16x^3+46x^2+42x+12)\\
                                         &=6x^2+10x+3>0
\end{align*}
for $x\geq 1$. Hence $B(k+1)>x+1$. By Claim~\ref{claim:Bmonotone}, for every $d\geq k+1$, $B(d)\geq B(k+1)>x+1$. Hence
$$
    \lfloor B(d)\rfloor+1\geq x+2>x+\frac{3}{2}=\frac{k+4}{4},
$$
contradicting \eqref{eq:b-upper}. Thus $k\equiv 2\pmod 4$ is impossible.

Now assume $k\equiv 0\pmod 4$.

Write $k=4x$ with $x\geq 3$. For $d=k+4=4x+4$, we have
$$
    \alpha(k+4)=\frac{(4x+4)(x+2)}{x(2x+5)},\qquad \alpha(k+4)+2=\frac{8x^2+22x+8}{x(2x+5)}.
$$
Therefore
$$
    B(k+4)=\frac{k+6}{\alpha(k+4)+2}=\frac{2x(2x+3)(2x+5)}{8x^2+22x+8}.
$$
We verify that $B(k+4)>x+1$. Since the denominator $8x^2+22x+8$ is positive, we have
\begin{align*}
    2x(2x+3)(2x+5)-(x+1)(8x^2+22x+8)&=(8x^3+32x^2+30x)-(8x^3+30x^2+30x+8)\\
                                    &=2x^2-8>0
\end{align*}
for $x\geq 3$. Hence $B(k+4)>x+1$. By Claim~\ref{claim:Bmonotone}, for every $d\geq k+4$, $B(d)\geq B(k+4)>x+1$. Hence
$$
    \lfloor B(d)\rfloor+1\geq x+2>x+1=\frac{k+4}{4},
$$
contradicting \eqref{eq:b-upper}. Thus for $k\equiv 0\pmod 4$ with $k\geq 12$, the only possible values of $d$ are
$$
    d=k+1,\ k+2,\ k+3.
$$
It remains to treat the cases $k=4$ and $k=8$. For these two values of $k$, the argument above does not restrict $d$ to $\{k+1,k+2,k+3\}$, so we must consider all $d\geq k+1$ directly.

\paragraph{Case 2.1.1.1. $k=4$.}\mbox{}\par

For $k=4$, we have
$$
    \frac{k+4}{4}=2,\qquad
    \alpha(d)=\frac{d(d-2)}{2(d-1)}.
$$
We first show that $r_1,r_m\geq 2$ for every $d\geq 5$. By \eqref{eq:rw-integer}, it suffices to prove $B(d)>1$. Indeed,
$$
    B(d)=\frac{d+2}{\alpha(d)+2}
    =\frac{2(d+2)(d-1)}{d(d-2)+4(d-1)}
    =\frac{2(d+2)(d-1)}{d^2+2d-4}.
$$
Since
$$
    2(d+2)(d-1)-(d^2+2d-4)
    =(2d^2+2d-4)-(d^2+2d-4)
    =d^2>0,
$$
we have $B(d)>1$. Hence $r_w\geq \lfloor B(d)\rfloor+1\geq 2$, so $r_1\geq 2$ and $r_m\geq 2$. Therefore $r_1+r_m\geq 4$. On the other hand, by \eqref{eq:r1rm-upper}, $r_1+r_m\leq \frac{k+4}{2}=\frac{8}{2}=4$. Thus equality holds: $r_1+r_m=4$ and $r_1=r_m=2$.

We now apply the $Q_1'$-technique. By \eqref{eq:r1primelower} and substituting $\alpha(d)=\frac{d(d-2)}{2(d-1)}$, we have
$$
    r_1'>\frac{d+3-2\cdot 2}{\alpha(d)}=\frac{2(d-1)^2}{d(d-2)}.
$$
Since
$$
    (d-1)^2-d(d-2)=(d^2-2d+1)-(d^2-2d)=1>0,
$$
we have
$$
    \frac{2(d-1)^2}{d(d-2)}>2
$$
for every $d\geq 5$. Therefore $r_1'>2$ and $r_1'+r_m>2+2=4$. By \eqref{eq:r1primermupper}, we also have $r_1'+r_m\leq \frac{k+4}{2}=4$, a contradiction. Hence $k=4$ is impossible.

\paragraph{Case 2.1.1.2. $k=8$.}\mbox{}\par

We have
$$
    \frac{k+4}{4}=3,\qquad
    \alpha(d)=\frac{2d(2d-8)}{8(2d-6)}=\frac{d(d-4)}{4(d-3)}.
$$
We need to consider all $d\geq k+1=9$. We split into three cases: $d=9,10,11$, $d=12$, and $d\geq 13$.

\begin{enumerate}[(i)]
\item $d=9,10,11$.

For $d=9$, we compute
$$
    \alpha(9)=\frac{9\cdot 5}{4\cdot 6}=\frac{15}{8}, \qquad
    B(9)=\frac{11}{\frac{15}{8}+2}=\frac{88}{31}\approx 2.839.
$$
Thus $\lfloor B(9)\rfloor +1=3$. Together with \eqref{eq:b-upper}, we get $\lfloor B(9)\rfloor+1=3=\frac{k+4}{4}$, so \eqref{eq:b-upper} holds with equality. Hence $r_1=r_m=3$.

For $d=10$, we compute
$$
    \alpha(10)=\frac{10\cdot 6}{4\cdot 7}=\frac{15}{7},\qquad
    B(10)=\frac{12}{\frac{15}{7}+2}=\frac{84}{29}\approx 2.897.
$$
Thus $\lfloor B(10)\rfloor +1=3$, and the same argument gives $r_1=r_m=3$.

For $d=11$, we compute
$$
    \alpha(11)=\frac{11\cdot 7}{4\cdot 8}=\frac{77}{32},\qquad
    B(11)=\frac{13}{\frac{77}{32}+2}=\frac{416}{141}\approx 2.950.
$$
Thus $\lfloor B(11)\rfloor +1=3$, and again $r_1=r_m=3$.

In all three cases, we apply the $Q_1'$-technique. Since $t_1=0$, By \eqref{eq:r1primelower} and substituting $\alpha(d)$, we have
$$
    r_1'> \frac{d+3-2\cdot 3}{\alpha(d)} = \frac{d-3}{\alpha(d)}.
$$
For $d=9$,
$$
    \alpha(9)=\frac{15}{8},\qquad
    \frac{d-3}{\alpha(9)}=\frac{6}{15/8}=\frac{16}{5}>3.
$$
For $d=10$,
$$
    \alpha(10)=\frac{15}{7},\qquad
    \frac{d-3}{\alpha(10)}=\frac{7}{15/7}=\frac{49}{15}>3.
$$
For $d=11$,
$$
    \alpha(11)=\frac{77}{32},\qquad
    \frac{d-3}{\alpha(11)}=\frac{8}{77/32}=\frac{256}{77}>3.
$$
Hence $r_1'>3$ in all three cases. By \eqref{eq:r1primermupper}, $r_1'+r_m\leq \frac{k+4}{2}=6$. But $r_1'>3$ and $r_m=3$, so $r_1'+r_m>3+3=6$, a contradiction.

\item $d=12$.

For $d=12$, we compute
$$
    \alpha(12)=\frac{12\cdot 8}{4\cdot 9}=\frac{8}{3},\qquad
    B(12)=\frac{14}{\frac{8}{3}+2}=3.
$$
Thus $\lfloor B(12)\rfloor+1=4>3=\frac{k+4}{4}$, contradicting
\eqref{eq:b-upper}. Hence $d=12$ is impossible.

\item $d\geq 13$.

For $d\geq 13$, we claim that $B(d)>3=\frac{k+4}{4}$. Indeed,
$$
    \alpha(d)=\frac{d(d-4)}{4(d-3)},\qquad
    B(d)=\frac{d+2}{\alpha(d)+2}.
$$
We compute
$$
    B(d)-3=\frac{d+2}{\alpha(d)+2}-3
    =\frac{d-4-3\alpha(d)}{\alpha(d)+2}.
$$
Since $\alpha(d)>0$, it suffices to show $d-4-3\alpha(d)>0$ for $d\geq 13$. Substituting $\alpha(d)=\frac{d(d-4)}{4(d-3)}$,
$$
    d-4-3\alpha(d)=d-4-\frac{3d(d-4)}{4(d-3)}
    =(d-4)\left(1-\frac{3d}{4(d-3)}\right)
    =(d-4)\cdot\frac{d-12}{4(d-3)}.
$$
For $d\geq 13$, both $d-4>0$ and $d-12>0$, so $d-4-3\alpha(d)>0$, hence $B(d)>3$. Therefore
$$
    \lfloor B(d)\rfloor+1\geq 4>3=\frac{k+4}{4},
$$
contradicting \eqref{eq:b-upper}. Thus $d\geq 13$ is impossible.

Combining the three subcases, we conclude that $k=8$ is impossible.

\end{enumerate}

\paragraph{Case 2.1.1.3. $k\geq 12$ and $d\in\{k+1,k+2,k+3\}$.}\mbox{}\par

Recall $k\equiv 0\pmod 4$ and write $k=4x$ with $x\geq 3$. We first show that $B(k+1)>x$. Indeed, for $d=k+1=4x+1$,
$$
    \alpha(k+1)=\frac{(4x+1)(2x+1)}{4x(x+1)},\qquad
    B(k+1)=\frac{4x(4x+3)(x+1)}{16x^2+14x+1}.
$$
Then
\begin{align*}
    B(k+1)-x
    =\frac{4x(4x+3)(x+1)-x(16x^2+14x+1)}{16x^2+14x+1}
    =\frac{14x+11}{16x^2+14x+1}>0.
\end{align*}
Hence $B(k+1)>x$. By Claim~\ref{claim:Bmonotone}, for every $d> k$, $B(d)\geq B(k+1)>x$. Therefore $\lfloor B(d)\rfloor+1\geq x+1$ for every $d> k$. By \eqref{eq:rw-integer}, $r_1,r_m\geq x+1$. By \eqref{eq:r1rm-upper}, $r_1+r_m\leq \frac{k+4}{2}=2x+2$. Together with $r_1,r_m\geq x+1$, we obtain $r_1=r_m=x+1$. We now apply the $Q_1'$-technique to $d=k+1,k+2,k+3$.

\begin{enumerate}[(i)]
\item $d=k+1$.

For $d=k+1=4x+1$, we have
$$
    \alpha(k+1)=\frac{(4x+1)(2x+1)}{4x(x+1)}.
$$
By \eqref{eq:r1primelower},
$$
    r_1'>\frac{d+3-2r_1}{\alpha(d)}
    =\frac{4x+1+3-2(x+1)}{\alpha(k+1)}
    =\frac{2x+2}{\alpha(k+1)}.
$$
Substituting $\alpha(k+1)=\frac{(4x+1)(2x+1)}{4x(x+1)}$,
$$
    r_1'>\frac{8x(x+1)^2}{(4x+1)(2x+1)}.
$$
We claim that $r_1'>x+1$. Indeed, this is equivalent to
$$
    \frac{8x(x+1)^2}{(4x+1)(2x+1)}>x+1
    \iff 8x(x+1)>(4x+1)(2x+1).
$$
The left-hand side is $8x^2+8x$, while the right-hand side is $8x^2+6x+1$. Their difference is $2x-1>0$ for $x\geq 1$. Hence $r_1'>x+1$. By \eqref{eq:r1primermupper}, $r_1'+r_m\leq \frac{k+4}{2}=2x+2$. But $r_1'>x+1$ and $r_m=x+1$, so $r_1'+r_m>2x+2$, a contradiction.

\item $d=k+2$.

For $d=k+2=4x+2$, we have
$$
    \alpha(k+2)=\frac{2(2x+1)(x+1)}{x(2x+3)}.
$$
By \eqref{eq:r1primelower},
$$
    r_1'>\frac{d+3-2r_1}{\alpha(d)}
    =\frac{4x+2+3-2(x+1)}{\alpha(k+2)}
    =\frac{2x+3}{\alpha(k+2)}.
$$
Substituting $\alpha(k+2)=\frac{2(2x+1)(x+1)}{x(2x+3)}$,
$$
    r_1'>\frac{x(2x+3)^2}{2(2x+1)(x+1)}.
$$
We claim that $r_1'>x+1$. Indeed, it suffices to prove
$$
    \frac{x(2x+3)^2}{2(2x+1)(x+1)}\geq x+1
    \iff x(2x+3)^2\geq 2(2x+1)(x+1)^2.
$$
The left-hand side is $4x^3+12x^2+9x$, while the right-hand side is $4x^3+10x^2+8x+2$. Their difference is $2x^2+x-2>0$ for $x\geq 1$. Hence $r_1'>x+1$. By \eqref{eq:r1primermupper}, $r_1'+r_m\leq \frac{k+4}{2}=2x+2$. But $r_1'>x+1$ and $r_m=x+1$, so $r_1'+r_m>2x+2$, a contradiction.

\item $d=k+3$.

For $d=k+3=4x+3$, we have
$$
    \alpha(k+3)=\frac{(4x+3)(2x+3)}{4x(x+2)}.
$$
By \eqref{eq:r1primelower},
$$
    r_1'>\frac{d+3-2r_1}{\alpha(d)}
    =\frac{4x+3+3-2(x+1)}{\alpha(k+3)}
    =\frac{2x+4}{\alpha(k+3)}.
$$
Substituting $\alpha(k+3)=\frac{(4x+3)(2x+3)}{4x(x+2)}$,
$$
    r_1'>\frac{8x(x+2)^2}{(4x+3)(2x+3)}.
$$
We claim that $r_1'>x+1$. Indeed, this is equivalent to
$$
    \frac{8x(x+2)^2}{(4x+3)(2x+3)}>x+1
    \iff 8x(x+2)^2>(x+1)(4x+3)(2x+3).
$$
The left-hand side is $8x^3+32x^2+32x$, while the right-hand side is $(x+1)(8x^2+18x+9)=8x^3+26x^2+27x+9$. Their difference is $6x^2+5x-9>0$ for $x\geq 1$. Hence $r_1'>x+1$. By
\eqref{eq:r1primermupper}, $r_1'+r_m\leq \frac{k+4}{2}=2x+2$. But $r_1'>x+1$ and $r_m=x+1$, so $r_1'+r_m>2x+2$, a contradiction.

Thus all three values $d=k+1,k+2,k+3$ are impossible. Therefore $k\equiv 0\pmod 4$ with $k\geq 12$ is impossible.
\end{enumerate}

\noindent\textbf{Case 2.1.2. $k$ is odd.}

Then
$$
    \alpha(d)=\frac{2d(2d-k+1)}{(k-1)(2d-k+3)}.
$$
Substituting $\alpha(d)$ into \eqref{eq:b-upper}, we analyze the following subcases.

\paragraph{Case 2.1.2.1. $k\equiv 3\pmod 4$.}\mbox{}\par

Write $k=4x+3$ with integer $x\geq 0$. We first compute $\alpha(k+1)$. For $d=k+1=4x+4$, we have
$$
    \alpha(k+1)=\frac{2(x+1)(2x+3)}{(2x+1)(x+2)},\qquad
    \alpha(k+1)+2=\frac{8x^2+20x+10}{(2x+1)(x+2)}.
$$
Consequently,
$$
    B(k+1)=\frac{k+3}{\alpha(k+1)+2}
    =\frac{4x+6}{\alpha(k+1)+2}
    =\frac{(2x+3)(2x+1)(x+2)}{4x^2+10x+5}.
$$
We now verify that $B(k+1)>x+1$. Since the denominator $8x^2+20x+10$ is positive, this is equivalent to
$$
    (2x+3)(2x+1)(x+2)>(x+1)(4x^2+10x+5).
$$
Expanding the left-hand side,
$$
    (2x+3)(2x+1)(x+2)=(4x^2+8x+3)(x+2)=4x^3+16x^2+19x+6,
$$
Expanding the right-hand side,
$$
    (x+1)(4x^2+10x+5)=4x^3+14x^2+15x+5.
$$
Their difference is $2x^2+4x+1>0$ for every $x\geq 0$. Hence $B(k+1)>x+1$.

By Claim~\ref{claim:Bmonotone}, for every $d\geq k+1$, $B(d)\geq B(k+1)>x+1$. Hence
$$
    \lfloor B(d)\rfloor+1\geq x+2>x+\frac{7}{4}=\frac{k+4}{4},
$$
contradicting \eqref{eq:b-upper}. Thus $k\equiv 3\pmod 4$ is impossible.

\paragraph{Case 2.1.2.2. $k\equiv 1\pmod 4$.}\mbox{}\par

Write $k=4x+1$ with integer $x\geq 1$. We first compute $B(d)$ for $d=k+1,k+2,k+3$.

\paragraph{Case 2.1.2.2.1. $d=k+1$.}\mbox{}\par

For $d=k+1=4x+2$, we have
$$
    \alpha(k+1)=\frac{2(2x+1)(x+1)}{x(2x+3)},\qquad
    B(k+1)=\frac{2x(x+1)(2x+3)}{4x^2+6x+1}.
$$
We first show that $B(k+1)>x$. Since $x>0$,
$$
    B(k+1)>x
    \iff 2(x+1)(2x+3)>4x^2+6x+1.
$$
The left-hand side is $4x^2+10x+6$, while the right-hand side is $4x^2+6x+1$. Their difference is $4x+5>0$, so $B(k+1)>x$. Hence $\lfloor B(k+1)\rfloor+1\geq x+1$. By \eqref{eq:rw-integer}, $r_1,r_m\geq x+1$. By \eqref{eq:r1rm-upper}, $r_1+r_m\leq \frac{k+4}{2}=2x+2.5$. Since $r_1+r_m$ is an integer, $r_1+r_m\leq 2x+2$. Together with $r_1,r_m\geq x+1$, we obtain $r_1=r_m=x+1$. We apply the $Q_1'$-technique. By \eqref{eq:r1primelower},
$$
    r_1'>\frac{d+3-2r_1}{\alpha(d)}
    =\frac{4x+2+3-2(x+1)}{\alpha(k+1)}
    =\frac{2x+3}{\alpha(k+1)}.
$$
Substituting $\alpha(k+1)=\frac{2(2x+1)(x+1)}{x(2x+3)}$,
$$
    r_1'>\frac{x(2x+3)^2}{2(2x+1)(x+1)}.
$$
We claim that $r_1'>x+1$. Indeed, this is equivalent to
$$
    x(2x+3)^2>2(2x+1)(x+1)^2.
$$
The left-hand side is $4x^3+12x^2+9x$, while the right-hand side is $4x^3+10x^2+8x+2$. Their difference is $2x^2+x-2>0$ for $x\geq 1$. Hence $r_1'>x+1$. By \eqref{eq:r1primermupper}, $r_1'+r_m\leq \frac{k+4}{2}=2x+2.5$. But $r_1'>x+1$ and $r_m=x+1$, so $r_1'+r_m>2x+2$. Since $r_1'+r_m$ is an integer, $r_1'+r_m\geq 2x+3>2x+2.5$, a contradiction. Hence $d=k+1$ is impossible.

\paragraph{Case 2.1.2.2.2. $d=k+2$.}\mbox{}\par

For $d=k+2=4x+3$, we have
$$
    \alpha(k+2)=\frac{(4x+3)(2x+3)}{4x(x+2)},\qquad \alpha(k+2)+2=\frac{16x^2+34x+9}{4x(x+2)}.
$$
Therefore
$$
    B(k+2)=\frac{k+4}{\alpha(k+2)+2}
    =\frac{4x(4x+5)(x+2)}{16x^2+34x+9}.
$$
We compare $B(k+2)$ with $x+1$. Since the denominator is positive,
$$
    B(k+2)>x+1
    \iff 4x(4x+5)(x+2)>(x+1)(16x^2+34x+9).
$$
Their difference is
\begin{align*}
    4x(4x+5)(x+2) - (x+1)(16x^2+34x+9)
    &=(16x^3+52x^2+40x)-(16x^3+50x^2+43x+9)\\
    &=2x^2-3x-9.
\end{align*}
\begin{enumerate}[(i)]
\item $x\geq 4$ (i.e., $k\geq 17$).

For $x\geq 4$, $2x^2-3x-9>0$. Hence $B(k+2)>x+1$. By Claim~\ref{claim:Bmonotone}, for every $d\geq k+2$, $B(d)\geq B(k+2)>x+1$. Hence $\lfloor B(d)\rfloor+1\geq x+2>x+\frac{5}{4}=\frac{k+4}{4}$, contradicting \eqref{eq:b-upper}. Hence $k\geq 17$ is impossible.

\item $x=1$.

We have $k=5$ and $d=7$. Then
$$
    \alpha(7)=\frac{2\cdot 7\cdot 10}{4\cdot 12}=\frac{35}{12},\qquad
    B(7)=\frac{9}{\alpha(7)+2}=\frac{9}{\frac{35}{12}+2}
    =\frac{108}{59}\approx 1.831.
$$
Thus $\lfloor B(7)\rfloor+1=2$. By \eqref{eq:r1rm-upper}, $r_1+r_m\leq \frac{k+4}{2}=\frac{9}{2}=4.5$. Since $r_1+r_m$ is an integer, $r_1+r_m\leq 4$. Also, by
\eqref{eq:rw-integer}, $r_1,r_m\geq 2$, so $r_1=r_m=2$. We apply the
$Q_1'$-technique. By \eqref{eq:r1primelower},
$$
    r_1'>\frac{d+3-2r_1}{\alpha(d)}
    =\frac{7+3-4}{\alpha(7)}
    =\frac{6}{35/12}=\frac{72}{35}>2.
$$
By \eqref{eq:r1primermupper}, $r_1'+r_m\leq \frac{k+4}{2}=4.5$. But $r_1'>2$ and $r_m=2$, so $r_1'+r_m>4$. Since $r_1'+r_m$ is an integer, $r_1'+r_m\geq 5>4.5$, a contradiction. Hence $k=5$, $d=7$ is impossible.

\item $x=2$.

We have $k=9$ and $d=11$. Then
$$
    \alpha(11)=\frac{2\cdot 11\cdot 14}{8\cdot 16}=\frac{77}{32},\qquad
    B(11)=\frac{13}{\alpha(11)+2}=\frac{13}{\frac{77}{32}+2}
    =\frac{416}{141}\approx 2.950.
$$
Thus $\lfloor B(11)\rfloor+1=3$. By
\eqref{eq:rw-integer}, $r_1,r_m\geq 3$. By \eqref{eq:r1rm-upper}, $r_1+r_m\leq \frac{k+4}{2}=\frac{13}{2}=6.5$.
Since $r_1+r_m$ is an integer, $r_1+r_m\leq 6$. So $r_1=r_m=3$. We apply the $Q_1'$-technique. By \eqref{eq:r1primelower},
$$
    r_1'>\frac{d+3-2r_1}{\alpha(d)}
    =\frac{11+3-6}{77/32}
    =\frac{256}{77}>3.
$$
By \eqref{eq:r1primermupper}, $r_1'+r_m\leq \frac{k+4}{2}=6.5$. But $r_1'>3$ and $r_m=3$, so $r_1'+r_m>6$. Since $r_1'+r_m$ is an integer, $r_1'+r_m\geq 7>6.5$, a contradiction. Hence $k=9$, $d=11$ is impossible.

\item $x=3$.

We have $k=13$ and $d=15$. Then
$$
    \alpha(15)=\frac{2\cdot 15\cdot 18}{12\cdot 20}=\frac{9}{4},\qquad B(15)=\frac{17}{\alpha(15)+2}=\frac{17}{\frac{9}{4}+2}=4.
$$
Thus $\lfloor B(15)\rfloor+1=5>4.25=\frac{k+4}{4}$, contradicting \eqref{eq:b-upper}. Hence $k=13$, $d=15$ is impossible.
\end{enumerate}

\paragraph{Case 2.1.2.2.3. $d\geq k+3$.}\mbox{}\par

For $d=k+3=4x+4$, we have
$$
    \alpha(k+3)=\frac{(4x+4)(x+2)}{x(2x+5)},\qquad \alpha(k+3)+2=\frac{8x^2+22x+8}{x(2x+5)}.
$$
Therefore
$$
    B(k+3)=\frac{k+5}{\alpha(k+3)+2}
    =\frac{x(4x+6)(2x+5)}{8x^2+22x+8}.
$$
We compare $B(k+3)$ with $x+1$. Since the denominator is positive,
$$
B(k+3)>x+1
    \iff x(4x+6)(2x+5)>(x+1)(8x^2+22x+8).
$$
Their difference is
$$
    x(4x+6)(2x+5)-(x+1)(8x^2+22x+8)=(8x^3+32x^2+30x)-(8x^3+30x^2+30x+8)=2x^2-8.
$$
\begin{enumerate}[(i)]
\item $x\geq 3$.

We have $k\geq 13$. For $x\geq 3$, $2x^2-8>0$. Hence $B(k+3)>x+1$. By Claim~\ref{claim:Bmonotone}, for every $d\geq k+3$, $B(d)\geq B(k+3)>x+1$. Hence $\lfloor B(d)\rfloor+1\geq x+2>x+\frac{5}{4}=\frac{k+4}{4}$, contradicting \eqref{eq:b-upper}. Hence $k\geq 13$ is impossible.

\item $x=1$.

We have $k=5$ and $d\geq 8$. Hence
$$
    \alpha(d)=\frac{2d(2d-4)}{4(2d-2)}=\frac{d(d-2)}{2(d-1)},\qquad \alpha(d)+2=\frac{d^2+2d-4}{2(d-1)}.
$$
Therefore,
$$
    B(d)=\frac{d+2}{\alpha(d)+2}
    =\frac{2(d+2)(d-1)}{d^2+2d-4}.
$$
We first show that $B(d)>1$ for every $d\geq 8$. Since
$$
    2(d+2)(d-1)-(d^2+2d-4)
    =2d^2+2d-4-(d^2+2d-4)
    =d^2>0,
$$
we have $B(d)>1$. Hence $\lfloor B(d)\rfloor+1\geq 2$ for every $d\geq 8$. By \eqref{eq:rw-integer}, $r_1,r_m\geq 2$. By \eqref{eq:r1rm-upper},
$$
    r_1+r_m\leq \frac{k+4}{2}=\frac{9}{2}=4.5.
$$
Since $r_1+r_m$ is an integer, $r_1+r_m\leq 4$. Together with $r_1,r_m\geq 2$, we obtain $r_1=r_m=2$. We apply the $Q_1'$-technique. By \eqref{eq:r1primelower},
$$
    r_1'>\frac{d+3-2r_1}{\alpha(d)}
    =\frac{d+3-4}{\alpha(d)}
    =\frac{d-1}{\alpha(d)}
    =\frac{2(d-1)^2}{d(d-2)}.
$$
We claim that $r_1'>2$. Indeed, this is equivalent to
$$
    (d-1)^2>d(d-2),
$$
and
$$
    (d-1)^2-d(d-2)=d^2-2d+1-(d^2-2d)=1>0.
$$
By \eqref{eq:r1primermupper}, $r_1'+r_m\leq \frac{k+4}{2}=4.5$. But $r_1'>2$ and $r_m=2$, so $r_1'+r_m>4$. Since $r_1'+r_m$ is an integer, $r_1'+r_m\geq 5>4.5$, a contradiction. Hence $k=5$, $d\geq 8$ is impossible.

\item $x=2$.

We have $k=9$, $d\geq 12$. Hence
$$
    \alpha(d)=\frac{2d(2d-8)}{8(2d-6)}=\frac{d(d-4)}{4(d-3)},\qquad \alpha(d)+2=\frac{d^2+4d-24}{4(d-3)}.
$$
Therefore,
$$
    B(d)=\frac{d+2}{\alpha(d)+2}
    =\frac{4(d+2)(d-3)}{d^2+4d-24}.
$$
We first show that $B(d)>3$ for every $d\geq 12$. We compute
$$
    B(d)-3=\frac{4(d+2)(d-3)-3(d^2+4d-24)}{d^2+4d-24}
    =\frac{d^2-10d+48}{d^2+4d-24}.
$$
For $d\geq 12$, the numerator is $d^2-10d+48=(d-5)^2+23>0$, so $B(d)>3$. Hence $\lfloor B(d)\rfloor+1\geq 4$ for every $d\geq 12$. Since $\frac{k+4}{4}=3.25$, we have $\lfloor B(d)\rfloor+1\geq 4>3.25=\frac{k+4}{4}$, contradicting \eqref{eq:b-upper}. Hence $k=9$, $d\geq 12$ is impossible.

\end{enumerate}

Thus all cases with $t_1 = t_m = 0$ lead to a contradiction.

\subsubsection{Case 2.2. $T=1$.}

Assume without loss of generality that $t_1=0$ and $t_m=1$ (otherwise reverse the path $P$). By \eqref{eq:main-ineq-general},
$$
    k+1> 2\cdot\frac{d+2-t_1}{\alpha(d)+2}
    +2\cdot\frac{d+2-t_m}{\alpha(d)+2}
    -3+T.
$$
Since $T=t_1+t_m=1$,
$$
    k+3>2\cdot\frac{2d+3}{\alpha(d)+2}.
$$
Multiplying both sides by $\alpha(d)+2>0$, $(k+3)(\alpha(d)+2)>4d+6$, which simplifies to
\begin{equation}\label{eq:T1-main}
(k+3)\alpha(d)>4d-2k.
\end{equation}

\paragraph{Case 2.2.1. $k$ is even.}\mbox{}\par

Then
$$
    \alpha(d)=\frac{2d(2d-k)}{k(2d-k+2)}.
$$
Substituting into
\eqref{eq:T1-main},
$$
    (k+3)\cdot\frac{2d(2d-k)}{k(2d-k+2)}>4d-2k.
$$
Multiplying both sides by $k(2d-k+2)>0$ and dividing by $2(2d-k)>0$, $d(k+3)>k(2d-k+2)$. So
\begin{equation}\label{eq:T1-even}
d(k-3)<k(k-2).
\end{equation}
For $k\geq 4$, this gives
$$
    d<\frac{k(k-2)}{k-3}=k+1+\frac{3}{k-3}.
$$
For even $k\geq 6$, one verifies that this forces $d=k+1$. The case $k=4$ is handled separately below.

\paragraph{Case 2.2.1.1. $k\geq 6$.}\mbox{}\par

By \eqref{eq:main-ineq-general}, $r_1+r_m\leq\frac{k+3}{2}$. Since $r_1+r_m$ is an integer and $\frac{k+3}{2}$ is a half-integer for even $k$, we have
$$
    r_1+r_m\leq\left\lfloor\frac{k+3}{2}\right\rfloor=\frac{k+2}{2}.
$$
On the other hand, for $d=k+1$, we have
$$
    \alpha(k+1)=\frac{2(k+1)(k+2)}{k(k+4)},\qquad \alpha(k+1)+2=\frac{4k^2+14k+4}{k(k+4)}.
$$
By \eqref{eq:rw-B},
$$
    r_1>\frac{(k+3)k(k+4)}{4k^2+14k+4},\qquad
    r_m>\frac{(k+2)k(k+4)}{4k^2+14k+4}.
$$
Hence
$$
    r_1+r_m>\frac{(2k+5)k(k+4)}{4k^2+14k+4}
    =\frac{2k^3+13k^2+20k}{4k^2+14k+4}.
$$
We compare this with $\frac{k+2}{2}$:
$$
    \frac{2k^3+13k^2+20k}{4k^2+14k+4}-\frac{k+2}{2}
    =\frac{4(k^2+2k-2)}{2(4k^2+14k+4)}>0
$$
for $k\geq 1$. Hence
$$
    r_1+r_m>\frac{k+2}{2},
$$
contradicting $r_1+r_m\leq\frac{k+2}{2}$. Therefore $k\geq 6$ even is impossible.

\paragraph{Case 2.2.1.2. $k=4$.}\mbox{}\par

For $k=4$, we have
$$
    \alpha(d)=\frac{d(d-2)}{2(d-1)},\qquad
    \alpha(d)+2=\frac{d^2+2d-4}{2(d-1)}.
$$
By \eqref{eq:rw-B}, for $w=1$ (recall $t_1=0$),
$$
    r_1>\frac{d+2}{\alpha(d)+2}
    =\frac{(d+2)\cdot 2(d-1)}{d^2+2d-4}.
$$
We claim that
$$
    \frac{(d+2)\cdot 2(d-1)}{d^2+2d-4}>\frac{d+2}{d+1}.
$$
Indeed, since $d+2>0$, this is equivalent to
$$
    2(d-1)(d+1)>d^2+2d-4.
$$
The left-hand side is $2d^2-2$, and the right-hand side is $d^2+2d-4$. Their difference is
$$
    2d^2-2-(d^2+2d-4)=d^2-2d+2=(d-1)^2+1>0.
$$
Hence $r_1>\frac{d+2}{d+1}>1$, so $r_1\geq 2$.

Similarly, for $w=m$ (recall $t_m=1$),
$$
    r_m>\frac{d+1}{\alpha(d)+2}
    =\frac{(d+1)\cdot 2(d-1)}{d^2+2d-4}.
$$
We claim that
$$
    \frac{(d+1)\cdot 2(d-1)}{d^2+2d-4}>1.
$$
Indeed, this is equivalent to
$$
    2(d+1)(d-1)>d^2+2d-4,
$$
which simplifies to
$$
    d^2-2d+2=(d-1)^2+1>0.
$$
Hence $r_m>1$, so $r_m\geq 2$.

Therefore $r_1+r_m\geq 4$. On the other hand, by
\eqref{eq:main-ineq-general},
$$
    r_1+r_m\leq\frac{k+4-T}{2}=\frac{4+4-1}{2}=3.5.
$$
But $r_1+r_m\geq 4>3.5$, a contradiction. Hence $k=4$ is impossible.

\paragraph{Case 2.2.2. $k$ is odd.}\mbox{}\par

Then
$$
    \alpha(d)=\frac{2d(2d-k+1)}{(k-1)(2d-k+3)}.
$$
Substituting into \eqref{eq:T1-main},
$$
    (k+3)\cdot\frac{2d(2d-k+1)}{(k-1)(2d-k+3)}>4d-2k.
$$
Multiplying both sides by $(k-1)(2d-k+3)>0$ and dividing by $2$,
$$
    d(k+3)(2d-k+1)>(2d-k)(k-1)(2d-k+3).
$$
Expanding and collecting like terms gives
\begin{equation}\label{eq:T1-odd}
    2(5-k)d^2+3(k-1)(k-3)d-k(k-1)(k-3)>0.
\end{equation}

\begin{enumerate}[(i)]
\item $k\geq 7$.

For $k\geq 7$, the coefficient $2(5-k)$ is negative, so the left-hand side of \eqref{eq:T1-odd} is a downward-opening parabola (or a linear function) in $d$. Its positive root is
$$
    d_0=\frac{-3(k-1)(k-3)+\sqrt{\Delta}}{4(5-k)},
$$
where
$$
    \Delta=9(k-1)^2(k-3)^2+8k(k-1)(k-3)(5-k).
$$
For general odd $k\geq 7$, we verify that $d_0<k+1$. Since $4(5-k)<0$, the inequality $d_0<k+1$ is equivalent to
$$
    \sqrt{\Delta}>3(k-1)(k-3)+4(5-k)(k+1).
$$
The right-hand side equals
$$
    3(k^2-4k+3)+4(-k^2+4k+5)=-k^2+4k+29.
$$
For $k\geq 9$, the right-hand side is negative, so the inequality holds trivially. For $k=7$, the right-hand side is $8$, and a direct computation gives $\sqrt{\Delta}=\sqrt{2496}\approx 49.96>8$. Hence $d_0<k+1$ for all odd $k\geq 7$. Since \eqref{eq:T1-odd} holds only for $d<d_0$, and $d\geq k+1>d_0$, inequality \eqref{eq:T1-odd} fails, a contradiction. Hence $k\geq 7$ odd is impossible.

\item $k=5$.

For $k=5$, we have
$$
    \alpha(d)=\frac{d(d-2)}{2(d-1)},\qquad
    \alpha(d)+2=\frac{d^2+2d-4}{2(d-1)}.
$$
By \eqref{eq:rw-B}, for $w=1$ (recall $t_1=0$) and $w=m$ (recall $t_m=1$),
$$
    r_1>\frac{d+2}{\alpha(d)+2},\qquad
    r_m>\frac{d+1}{\alpha(d)+2}.
$$
Since $d+2>d+1$, we have $r_1>r_m$. It therefore suffices to show $r_m>1$. Indeed,
$$
    r_m>\frac{d+1}{\alpha(d)+2}
    =\frac{2(d+1)(d-1)}{d^2+2d-4}.
$$
We claim that
$$
    \frac{2(d+1)(d-1)}{d^2+2d-4}>1.
$$
This is equivalent to
$$
    2(d+1)(d-1)>d^2+2d-4,
$$
which simplifies to $d^2-2d+2=(d-1)^2+1>0$. Hence $r_m>1$, and therefore $r_1>r_m>1$. So $r_1\geq 2$ and $r_m\geq 2$. By \eqref{eq:main-ineq-general}, $r_1+r_m\leq\frac{5+4-1}{2}=4$. Hence $r_1=r_m=2$. We apply the $Q_1'$-technique. Since $t_1=0$, we have $S_1\subseteq R_1^-$. By \eqref{eq:r1primelower},
$$
    r_1'>\frac{d+3-2r_1}{\alpha(d)}=\frac{d-1}{\alpha(d)}
    =\frac{2(d-1)^2}{d(d-2)}.
$$
For $d\geq 6$, we have $\frac{2(d-1)^2}{d(d-2)}>2$. Hence $r_1'>2$. By
\eqref{eq:r1primermupper}, $r_1'+r_m\leq\frac{5+4-1}{2}=4$.
But $r_1'>2$ and $r_m=2$, so $r_1'+r_m>4$, a contradiction. Hence $k=5$ is impossible.

\item $k=3$.

For $k=3$, we have $\alpha(d)=d-1$. By \eqref{eq:rw-B},
$$
    r_1>\frac{d+2}{d+1}>1\implies r_1\geq 2,\qquad
    r_m>\frac{d+1}{d+1}=1\implies r_m\geq 2.
$$
By \eqref{eq:main-ineq-general}, $r_1+r_m\leq\frac{3+4-1}{2}=3$. But $r_1+r_m\geq 4>3$, a contradiction. Hence $k=3$ is impossible.
\end{enumerate}

Therefore all cases with $T=1$ lead to a contradiction.

\subsubsection{Case 2.3. \texorpdfstring{$T \geq 2$}{T >= 2}.}

By \eqref{eq:main-ineq-general},
$$
    k+1> 2\cdot\frac{d+2-t_1}{\alpha(d)+2}+2\cdot\frac{d+2-t_m}{\alpha(d)+2}-3+T.
$$
Rearranging,
\begin{equation}\label{eq:T-general}
    \frac{2(2d+4-T)}{\alpha(d)+2}<k+4-T.
\end{equation}
Since $T\geq 2$, we have
$$
    2d+4-T\geq 2d+6-2T=2(d+3-T).
$$
Combining this with \eqref{eq:T-general},
$$
    \frac{4(d+3-T)}{\alpha(d)+2}
<k+4-T.
$$
Multiplying both sides by $\alpha(d)+2>0$,
$$
    4(d+3-T)<(k+4-T)(\alpha(d)+2).
$$
Expanding and rearranging, $$4d+12-4T<(k+4)(\alpha(d)+2)-T(\alpha(d)+2).$$
So
$$
    T(\alpha(d)-2)<(k+4)(\alpha(d)+2)-4d-12.
$$
Since $T\geq 2$, $2(\alpha(d)-2)<(k+4)(\alpha(d)+2)-4d-12,$ which simplifies to
\begin{equation}\label{eq:T2-main}
    (k+2)\alpha(d)>4d-2k.
\end{equation}
We now substitute the definition of $\alpha(d)$ and derive a contradiction in both parity cases.

\paragraph{Case 2.3.1. $k$ is even.}\mbox{}\par

Then
$$
    \alpha(d)=\frac{2d(2d-k)}{k(2d-k+2)}.
$$
Substituting into \eqref{eq:T2-main},
$$
    (k+2)\cdot\frac{2d(2d-k)}{k(2d-k+2)}>4d-2k.
$$
Multiplying both sides by $k(2d-k+2)>0$,
$$
    2d(k+2)(2d-k)>k(2d-k+2)(4d-2k).
$$
Since $2d-k>0$ and $4d-2k=2(2d-k)$, we may divide both sides by $2(2d-k)>0$ to obtain
$$
    d(k+2)>k(2d-k+2).
$$
Expanding,
$$
    dk+2d>2dk-k^2+2k.
$$
So
$$
    0>dk-2d-k^2+2k=(k-2)(d-k).
$$
Since $k\geq 4$ and $d>k$, we have $(k-2)(d-k)>0$, a contradiction. Thus $k$ even is impossible.

\paragraph{Case 2.3.2. $k$ is odd.}\mbox{}\par

Then
$$
    \alpha(d)=\frac{2d(2d-k+1)}{(k-1)(2d-k+3)}.
$$
Substituting into \eqref{eq:T2-main},
$$
    (k+2)\cdot\frac{2d(2d-k+1)}{(k-1)(2d-k+3)}>4d-2k.
$$
Multiplying both sides by $(k-1)(2d-k+3)>0$, and dividing both sides by $2$,
$$
    d(k+2)(2d-k+1)>(2d-k)(k-1)(2d-k+3).
$$
Subtracting the right-hand side from the left-hand side and collecting like terms gives
\begin{equation}\label{eq:T2-odd}
    2(4-k)d^2+(k-1)(3k-8)d-k(k-1)(k-3)>0.
\end{equation}

\paragraph{Case 2.3.2.1. $k\geq 5$.}\mbox{}\par

For $k\geq 5$, the coefficient $2(4-k)$ is negative, so the left-hand side of \eqref{eq:T2-odd} is a downward-opening parabola in $d$. Its positive root is
$$
    d_0=\frac{-(k-1)(3k-8)+\sqrt{\Delta}}{4(4-k)},
$$
where
$$
    \Delta=(k-1)^2(3k-8)^2+8k(k-1)(k-3)(4-k).
$$
For $k=5$ and $k=7$, a direct computation gives $d_0<k+1$:
$$
    k=5:\ d_0\approx 1.615<6=k+1,\qquad
    k=7:\ d_0\approx 2.725<8=k+1.
$$
For general odd $k\geq 5$, we claim that
$$
    d_0<\frac{k-1}{2}<k+1.
$$
To see this, note that $4(4-k)<0$, so $d_0<\frac{k-1}{2}$ is equivalent to
$$
    \sqrt{\Delta}>(k-1)(3k-8)+2(4-k)(k-1)=k(k-1).
$$
Squaring both sides, it suffices to prove
$$
    \Delta>k^2(k-1)^2.
$$
Computing the difference,
\begin{align*}
    \Delta-k^2(k-1)^2
    &=(k-1)^2\bigl[(3k-8)^2-k^2\bigr]+8k(k-1)(k-3)(4-k)\\
    &=8(k-1)^2(k-2)(k-4)+8k(k-1)(k-3)(4-k)\\
    &=8(k-1)\bigl[(k-1)(k-2)(k-4)+k(k-3)(4-k)\bigr].
\end{align*}
The bracket equals
$$
    (k^3-7k^2+14k-8)+(-k^3+7k^2-12k)=2(k-4)>0
$$
for $k\geq 5$. Hence $\Delta>k^2(k-1)^2$, so $\sqrt{\Delta}>k(k-1)$ and therefore
$$
    d_0<\frac{k-1}{2}<k+1.
$$
Since \eqref{eq:T2-odd} holds only for $d<d_0$, and $d\geq k+1>d_0$, inequality \eqref{eq:T2-odd} fails, a contradiction.

\paragraph{Case 2.3.2.2. $k=3$.}\mbox{}\par

For $k=3$, we have $\alpha(d)=d-1$. By \eqref{eq:main-ineq-general},
$$
    7-T>\frac{2(2d+4-T)}{d+1}.
$$
Multiplying both sides by $d+1>0$, $(7-T)(d+1)>2(2d+4-T)$, which simplifies to
$$
    T(d-1)<3d-1.
$$
So
$$
    T<\frac{3d-1}{d-1}=3+\frac{2}{d-1}.
$$
Since $d\geq 4$, we have $\frac{2}{d-1}\leq\frac{2}{3}<1$, hence $T<4$. Together with $T\geq 2$, this gives $T\in\{2,3\}$.

Assume without loss of generality that $t_1\leq t_m$. Since $T\in\{2,3\}$, we have $t_1\leq 1$. By \eqref{eq:rw-B},
$$
    r_1>\frac{d+2-t_1}{d+1}\geq\frac{d+1}{d+1}=1.
$$
So $r_1\geq 2$. Similarly, for $w=m$, we have $t_m\leq 3$. By \eqref{eq:rw-B},
$$
    r_m>\frac{d+2-t_m}{d+1}\geq\frac{d-1}{d+1}\geq\frac{3}{5}.
$$
So $r_m\geq 1$. Therefore $r_1+r_m\geq 3$. On the other hand, by \eqref{eq:main-ineq-general},
$$
    r_1+r_m\leq\frac{7-T}{2}.
$$
For $T=2$, this gives $r_1+r_m\leq 2.5$, hence $r_1+r_m\leq 2$. For $T=3$, this gives $r_1+r_m\leq 2$. In both cases, we obtain $r_1+r_m\leq 2$, contradicting $r_1+r_m\geq 3$. Hence $k=3$, $T\geq 2$ is impossible.

Thus $k$ odd is impossible. Therefore all cases with $t_1+t_m\geq 2$ lead to a contradiction. This completes the proof of Theorem~\ref{theorem:main}.
\end{proof}

\bibliographystyle{amsplain}
\bibliography{refs}
\end{document}